\documentstyle[amsfonts, amssymb, latexsym, epsfig, epic, eepic ,12pt]{amsart}

\newtheorem{theorem}{Theorem}[section]
\newtheorem{lemma}[theorem]{Lemma}
\newtheorem{proposition}[theorem]{Proposition}
\newtheorem{corollary}[theorem]{Corollary}
\newtheorem{definition}[theorem]{Definition}

\def\C{{\mbox{\rm\kern.24em
\vrule width.03em height1.43ex depth-.052ex \kern-.26em C}}}
\def\QSet{\mbox{\rm\kern.24em
\vrule width.03em height1.48ex depth-.051ex \kern-.26em Q}}
\def\Z{{\bf Z}}
\def\R{{\mbox{\rm I\kern-.22em R}}}
\def\P{{\bf P}}
\def\Q{{\bf Q}}
\def\T{{\bf T}}
\def\D{{\bf D}}
\def\supp{{\rm supp}}
\def\size{{\rm size}}
\def\diam{{\rm diam}}
\def\energy{{\rm energy}}
\def\modenergy{\widetilde{\rm energy}}

\def\pv{{\vec P}}
\def\qv{{\vec Q}}
\def\Pv{{\vec{\bf P}}}
\def\Qv{{\vec{\bf Q}}}
\def\dist{{\rm dist}}

\def\111{\gamma}

\def\be#1{\begin{equation}\label{#1}}
\def\bas{\begin{align*}}
\def\eas{\end{align*}}
\def\bi{\begin{itemize}}
\def\ei{\end{itemize}}
\newenvironment{proof}{\noindent {\bf Proof:} }{\endprf\par}
\def \endprf{\hfill  {\vrule height6pt width6pt depth0pt}\medskip}
\def\emph#1{{\it #1}}

\title{The Bi-Carleson operator}

\author{Camil Muscalu}
\address{Department of Mathematics, Cornell University, Ithaca, NY 14853}
\email{camil@@math.cornell.edu}

\address{Current Address: School of Mathematics, IAS Princeton, NJ 08540}
\email{camil@@math.ias.edu}

\author{Terence Tao}
\address{Department of Mathematics, UCLA, Los Angeles CA 90095-1555 }
\email{tao@@math.ucla.edu}

\author{Christoph Thiele}
\address{Department of Mathematics, UCLA, Los Angeles CA 90095-1555}
\email{thiele@@math.ucla.edu}

\include{psfig}
\begin{document}

\begin{abstract}  
We prove $L^p$ estimates (Theorem \ref{main}) for the Bi-Carleson operator defined below.  The methods used 
are essentially based on the treatment of the Walsh analogue of the operator in the 
prequel \cite{mtt:walshbicarleson} of this paper, but with additional technicalities due to the fact that in the Fourier model one cannot obtain perfect 
localization in both space and frequency.
\end{abstract}

\maketitle

\section{introduction}

The maximal Carleson operator is the sub-linear operator defined by

\begin{equation}
C(f)(x):= \sup_N\left| \int_{\xi<N}\widehat{f}(\xi) e^{2\pi i x \xi} d\xi \right|,
\end{equation}
where $f$ is a Schwartz function on the real line $\R$ and the Fourier transform is defined by

$$\widehat{f}(\xi):= \int_{\R} f(x) e^{-2 \pi i x \xi} dx .$$
The following statement of Carleson and Hunt \cite{carleson}, \cite{hunt} is a classical theorem in Fourier analysis:

\begin{theorem}\label{carleson-fourier}
The operator $C$ maps $L^p\rightarrow L^p$, for every $1<p<\infty$.
\end{theorem}
This result, in the particular weak type $(2,2)$ special case, was the main ingredient in the proof of Carleson's
famous theorem which states that the Fourier series of a function in $L^2(\R/\Z)$ converges pointwise almost everywhere.

The bilinear Hilbert transform is an operator which can be essentially written as

\begin{equation}
B(f_1, f_2)(x):= \int_{\xi_1 < \xi_2} \widehat{f_1}(\xi_1) \widehat{f_2}(\xi_2) e^{2\pi i x (\xi_1 + \xi_2)} d\xi_1 d\xi_2,
\end{equation}
where $f_1, f_2$ are test functions on $\R$. From the work of Lacey and Thiele  \cite{laceyt} we have the following $L^p$ estimates on $B$:

\begin{theorem}\label{bht} 
$B$ maps $L^{p_1} \times L^{p_2} \to L^{p'_3}$ whenever $1 < p_1, p_2 \leq \infty$, $1/p_1 + 1/p_2 = 1/p'_3$, and $2/3 < p'_3 < \infty$.
\end{theorem}
The purpose of the present paper is to study the $L^p$ boundedness properties of the Bi-Carleson operator defined by

\begin{equation}
T(f_1, f_2)(x):= \sup_N \left|\int_{\xi_1 < \xi_2 < N} \widehat{f_1}(\xi_1) \widehat{f_2}(\xi_2) e^{2\pi i x (\xi_1 + \xi_2)} d\xi_1 d\xi_2
\right|.
\end{equation}
Our main theorem is the following:

\begin{theorem}\label{main}
The Bi-Carleson operator $T$ defined above maps $L^{p_1} \times L^{p_2} \to L^{p'_3}$ as long as $1 < p_1, p_2 \leq \infty$, $1/p_1 + 1/p_2 = 1/p'_3$, and $2/3 < p'_3 < \infty$.
\end{theorem}

%\begin{figure}[htbp]\centering
%\psfig{figure=name2.eps, height=5in,width=5.6in}
%\caption{The polytope $[A_1 \ldots A_{12}]$}
%\label{fig}
%\end{figure}

This Bi-Carleson operator can be thought of as being a hybrid of the Carleson operator and the bilinear Hilbert transform. In fact, our main theorem above
implies Carleson-Hunt's theorem and Lacey-Thiele's theorem as special cases (in Section 11 we also give an expository proof of Carleson-Hunt theorem).

The interesting and beautiful fact about this operator (and its ``relative'' bi-est operator studied in \cite{mtt:walshbiest}, \cite{mtt:fourierbiest}) is that,
unlike the Carleson operator or the bilinear Hilbert transform, it has a very special $\it{biparameter}$ 
structure not seen among the previously studied operators in harmonic 
analysis. As the careful reader will notice, understanding this structure is the main challenge of the paper.

In \cite{mtt:walshbicarleson} Theorem \ref{main} was proved for a Walsh-Fourier analogue $T_{walsh,\P,\Q}$ of $T$.  From the point of view of time-frequency analysis the two operators are closely related, however the Walsh model is easier to analyze technically because it is possible to localize perfectly in both space and frequency simultaneously.  In the Fourier case one has to deal with several ``Schwartz tails'' which introduce additional difficulties.  For instance, in the Walsh model an inner product $\langle \phi_P, \phi_Q \rangle$ of wave packets vanishes unless the spatial intervals $I_P$ and $I_Q$ are nested; however in the Fourier model one needs to consider the case when $I_P$ and $I_Q$ are separated (although the estimates improve rapidly with the relative separation of $I_P$ and $I_Q$).

The study of $T$ clearly reduces to the study of its linearized version defined by

\begin{equation}\label{linearT}
T(f_1, f_2)(x):= \int_{\xi_1 < \xi_2 < N(x)} \widehat{f_1}(\xi_1) \widehat{f_2}(\xi_2) e^{2\pi i x (\xi_1 + \xi_2)} d\xi_1 d\xi_2,
\end{equation}
where ''$x\rightarrow N(x)$'' is an arbitrary function.

We recall now a theorem from \cite{lim} which will take care of an error term later on.

Let $m$ be a bounded function on $\R^2$, satisfying the classical Marcinkiewicz-Mihlin-H\"{o}rmander condition

\begin{equation}\label{mmh}
|\partial_{\vec{\xi}}^\alpha m(\vec{\xi})| 
\lesssim |\vec{\xi}|^{-|\alpha|}
\end{equation}
for sufficiently many multiindices $\alpha$ and define the multiplier operator
\[
T_m(f_1,f_2)(x):=\int_{\R^2} m(\xi_1,\xi_2)
\widehat{f}_1(\xi_1)\widehat{f}_2(\xi_2)
e^{2\pi ix(\xi_1+\xi_2)}\,d\xi_1 d\xi_2.
\]

Similarly, we consider the Carleson type operator associated to the symbol $m$ and defined by

\begin{equation}\label{carlesonm}
C_m(f_1, f_2)(x):= \sup_{\vec{N}\in\R^2} \left|T_{\tau_{\vec{N}}m}(f_1, f_2)(x)\right|,
\end{equation}
where $\tau_{\vec{N}}m(\vec{\xi}):= m(\vec{\xi}-\vec{N})$.

In \cite{lim} the following theorem has been proven:

\begin{theorem}\label{xc}
The operator $T_m$ maps $L^{p_1} \times L^{p_2} \to L^{p'_3}$ as long as $1 < p_1, p_2 \leq \infty$, $1/p_1 + 1/p_2 = 1/p'_3$, and $0 < p'_3 < \infty$.
\end{theorem}

As in \cite{mtt:walshbicarleson} it shall be convenient to split our linearized operator $T$ into two pieces plus an ``error term'', which is an operator of the 
form (\ref{carlesonm}).
Specifically, fix $N\in\R$ and construct the symbols $m'_N, m''_N, m'''_N$ such that
$$\supp (m'_N)\subseteq \{\xi_1<\xi_2; \frac{\xi_1+\xi_2}{2}<N\},$$ $m'_N$ is smooth away from the line $\xi_1=\xi_2$ and equals $1$ on $\xi_1=\xi_2$,
$$ \supp (m''_N)\subseteq \{\xi_1<\xi_2<N\},$$ $m''_N$ is smooth away from the line $\xi_2=N$ and equals $1$ on $\xi_2=N$, while
$$m'''_N:=\chi_{\xi_1<\xi_2<N} - m'_N - m''_N$$ is smooth away from the vertex $(N,N)$ of the cone $\{\xi_1<\xi_2<N\}$ and moreover, 
$\tau_{-(N,N)}m'''_N$ satisfies (\ref{mmh}) (with the corresponding constants independent of $N$).

The multipliers $m'_N$, $m''_N$ and $m'''_N$ will be carefully constructed in the forthcoming sections. Consequently, our operator
$T$ in (\ref{linearT}) can be written as

$$T = T' + T'' + T'''.$$
Since $T'''$ can be estimated by Theorem \ref{xc}, it is enough to prove that both $T'$ and $T''$ satisfy the conclusion of the main Theorem (\ref{main}).

While the current article is mostly self contained, a knowledge of the prequel \cite{mtt:walshbicarleson} should help. However, we believe that the overall presentation
of the ideas related to this problem, has been improved in the meantime. In particular, we managed to
avoid the use of the complicated mixed sizes ``$\size (f,g)$'' and energies ``$\energy(f,g)$'' which appeared in the second half of \cite{mtt:walshbicarleson}.

The first author was partially supported by  NSF.
The second author is a Clay Prize Fellow and is supported by a grant from
Packard Foundations. The third author was partially supported
by a Sloan Fellowship and by NSF.

\section{Notation}\label{notation-sec}

In this section we set out some general notation used throughout the paper.  

We use $A\lesssim B$ to denote the statement that
$A\leq CB$ for some large constant $C$, and $A\ll B$ to denote the
statement that $A\leq C^{-1}B$ for some large constant $C$. We will sometime write $A\sim B$ and this means that $A\lesssim B\lesssim A.$ Given any interval (or cube) $I$, we let $|I|$ denote the measure of $I$, and $cI$ denotes the interval (or cube) with the same center as $I$ but $c$ times the side-length.

Given a spatial interval $I$, we shall define the approximate cutoff function $\tilde \chi_I$ by
$$ \tilde \chi_I(x) := (1 + (\frac{|x-x_I|}{|I|})^2)^{-1/2},$$
where $x_I$ is the center of $I$.

A collection $\{\omega\}$ of intervals is said to be \emph{lacunary around the frequency $\xi$} if we have  $\dist(\xi, \omega) \sim |\omega|$ for all $\omega$ in the 
collection.

\section{interpolation}\label{interp-sec}

In this section we review the interpolation theory from \cite{mtt:multilinear} which 
allows us to reduce multi-linear $L^p$ estimates such as those in 
Theorem \ref{main} to certain ``restricted  type'' estimates.

To prove the $L^p$ estimates on $T$ it is convenient to use 
duality and introduce the trilinear form $\Lambda$ associated to $T$ via the formula
$$
\Lambda(f_1, f_2,f_3) := \int_{\R} 
T(f_1, f_2)(x)f_3(x) dx.
$$
Similarly, define $\Lambda'$ and $\Lambda''$ associated with $T'$ and $T''$ respectively.  The statement that $T$ is bounded from 
$L^{p_1} \times L^{p_2}$ to $L^{p'_3}$ is then equivalent to 
$\Lambda$ being bounded on 
$L^{p_1} \times L^{p_2} \times L^{p_3}$ if 
$1 \le  p_3' < \infty$.  For $p_3'<1$ this simple duality relationship breaks down, however the interpolation arguments in \cite{mtt:multilinear} 
will allow us to reduce our desired estimate to certain ``restricted
type'' estimates on $\Lambda$. As in \cite{mtt:walshbicarleson} we find more
convenient to work with the quantities $\alpha_i=1/p_i$, $i=1,2,3$,
where $p_i$ stands for the exponent of $L^{p_i}$.

\begin{definition}
A tuple $\alpha=(\alpha_1,\alpha_2,\alpha_3)$ is called admissible, if

\[-\infty <\alpha_i <1\]
for all $1\leq i\leq 3$,

\[\sum_{i=1}^3\alpha_i=1\]
and there is at most one index $j$ such that $\alpha_j<0$. We call
an index $i$ good if $\alpha_i\geq 0$, and we call it bad if
$\alpha_i<0$. A good tuple is an admissible tuple without bad index, a
bad tuple is an admissible tuple with a bad index.
\end{definition}

\begin{definition}
Let $E$, $E'$ be sets of finite measure. We say that $E'$ is a major 
subset of $E$ if $E'\subseteq E$ and $|E'|\geq\frac{1}{2}|E|$.
\end{definition} 

\begin{definition}
If $E$ is a set of finite measure, we denote by $X(E)$ the space of
all measurable complex-valued functions $f$ supported on $E$ and such that $\|f\|_{\infty}\leq
1$.
\end{definition}

\begin{definition}
If $\alpha=(\alpha_1,\alpha_2,\alpha_3)$ is an admissible bad tuple
with bad index $j$, we say
that a $3$-linear form $\Lambda$ is of restricted type $\alpha$
if for every sequence $E_1, E_2, E_3$ of subsets of $\R$ with
finite measure, there exists a major subset $E'_j$ of $E_j$ such that

\[|\Lambda(f_1, f_2, f_3)|\lesssim |E|^{\alpha}\]
for all functions $f_i\in X(E'_i)$, $i=1,2,3$, where we adopt the
convention $E'_i =E_i$ for good indices $i$, and $|E|^{\alpha}$
is a shorthand for

\[|E|^{\alpha}=|E_1|^{\alpha_1}|E_2|^{\alpha_2}
|E_3|^{\alpha_3}.\]
If $\alpha=(\alpha_1,\alpha_2,\alpha_3)$ is an admissible good tuple, 
we say
that a $3$-linear form $\Lambda$ is of restricted type $\alpha$ if
there exists $j$ such that
for every sequence $E_1, E_2, E_3$ of subsets of $\R$ with
finite measure, there exists a major subset $E'_j$ of $E_j$ such that

\[|\Lambda(f_1, f_2, f_3)|\lesssim |E|^{\alpha}\]
for all functions $f_i\in X(E'_i)$, $i=1,2,3$, where this time we adopt the
convention $E'_i =E_i$ for the indices $i\neq j$.

\end{definition}

Let us consider now the $2$-dimensional affine hyperspace
\[
S:=\{(\alpha_1,\alpha_2,\alpha_3)\in\R^3\,
|\,\alpha_1 + \alpha_2 + \alpha_3=1\}.
\]

%\begin{figure}[htbp]\centering
%\psfig{figure=3.eps,height=4.4in,width=5in}
%\caption{The polytope $[A_1...A_{6}]$}
%\label{fig}
%\end{figure}

The points $A_1,...,A_6$ belong to $S$ and have the following coordinates:

\[
\begin{array}{llll}
A_1:(-\frac 12,1,\frac 12)  &  A_2:(\frac 12,1,-\frac 12)  &  
A_3:(1,\frac 12,-\frac 12)   \\  
\ &\ &\ & \ \\
A_4:(1,-\frac 12,\frac 12)  &  A_5:(\frac 12,-\frac 12,1) & 
A_6:(-\frac 12,\frac 12,1).   \\  
\end{array}
\]
The points $M_{12}, M_{34}, M_{56}$ are midpoints of their
corresponding segments and have the coordinates
$M_{12}:(0,1,0)$, $M_{34}:(1,0,0)$, $M_{56}:(0,0,1)$. Also, the point $A$ has the coordinates $(1,1,-1)$.

\setlength{\unitlength}{1.1mm}
\begin{figure}\label{hexagon}
\caption{Hexagon}
\begin{picture}(72,78)
%\put(46,40){\circle*{0.8}}
\put(46,60){\circle*{0.8}}
\put(28,30){\circle*{0.8}}
\put(64,30){\circle*{0.8}}
%\put(46,30){\circle*{1.6}}
\put(19,45){\circle*{0.8}}
\put(28,60){\circle*{0.8}}
\put(64,60){\circle*{0.8}}
\put(73,45){\circle*{0.8}}
\put(37,15){\circle*{0.8}}
\put(55,15){\circle*{0.8}}
\put(81.5,60){\circle*{0.8}}
%\put(55,45){\circle*{1.6}}
\put(27,62){$A_1$}
\put(45,62){$M_{12}$}
\put(63,62){$A_2$}
\put(80.5,62){$A$}
\put(13,44){$A_6$}
\put(74,44){$A_3$}
\put(20,29){$M_{56}$}
\put(65,29){$M_{34}$}
\put(31,14){$A_5$}
\put(56,14){$A_4$}

%\put(5,28){$(1,0,0)$}
%\put(73,28){$(0,1,0)$}
%\put(40,64){$(0,0,1)$}

%\put(27,50){$a$}
%\put(63,50){$b$}
%\put(45,35){$c$}
%\put(45,20){$d$}
%\put(45,7){$e$}

\put(10,60){\line(1,0){71.4}}
%\multiput(10,60)(6,0){2}{\line(1,0){4}}
%\multiput(29,60)(6,0){6}{\line(1,0){4}}
%\multiput(82,60)(-6,0){2}{\line(-1,0){4}}
%\put(19.4,45){\line(1,0){53.2}}
\put(28.4,30){\line(1,0){35.2}}
\put(38,15){\line(1,0){17.5}}
\put(10,60){\line(3,-5){35.8}}
\put(19.,45){\line(3,5){8.6}}
%\multiput(19.5,44.17)(3,-5){6}{\line(3,-5){2}}
%\multiput(46,0)(-3,5){2}{\line(-3,5){2}}
%\put(28.2,59.67){\line(3,-5){26.6}}
\put(46.2,59.67){\line(3,-5){17.6}}
\put(73,45){\line(-3,5){8.6}}
%\multiput(64.5,59.17)(3,-5){3}{\line(3,-5){2}}
%\multiput(19.5,45.83)(3,5){3}{\line(3,5){2}}
\put(28.2,30.33){\line(3,5){17.6}}
%\put(37.2,15.33){\line(3,5){26.6}}
\put(46,0){\line(3,5){35.8}}
%\multiput(46,0)(3,5){2}{\line(3,5){2}}
%\multiput(55.5,15.83)(3,5){6}{\line(3,5){2}}
%\multiput(82,60)(-3,-5){2}{\line(-3,-5){2}}

%\put(21,45){\line(3,5){8}}
%\put(23,45){\line(3,5){7}}
%\put(25,45){\line(3,5){6}}
%\put(27,45){\line(3,5){5}}
%\put(29,45){\line(3,5){4}}
%\put(31,45){\line(3,5){3}}
%\put(33,45){\line(3,5){2}}
%\put(35,45){\line(3,5){1}}

%\put(39,45){\line(3,-5){8}}
%\put(41,45){\line(3,-5){7}}
%\put(43,45){\line(3,-5){6}}
%\put(45,45){\line(3,-5){5}}
%\put(47,45){\line(3,-5){4}}
%\put(49,45){\line(3,-5){3}}
%\put(51,45){\line(3,-5){2}}
%\put(53,45){\line(3,-5){1}}

%\put(57,45){\line(3,5){8}}
%\put(59,45){\line(3,5){7}}
%\put(61,45){\line(3,5){6}}
%\put(63,45){\line(3,5){5}}
%\put(65,45){\line(3,5){4}}
%\put(67,45){\line(3,5){3}}
%\put(69,45){\line(3,5){2}}
%\put(71,45){\line(3,5){1}}

%\put(39,15){\line(3,5){8}}
%\put(41,15){\line(3,5){7}}
%\put(43,15){\line(3,5){6}}
%\put(45,15){\line(3,5){5}}
%\put(47,15){\line(3,5){4}}
%\put(49,15){\line(3,5){3}}
%\put(51,15){\line(3,5){2}}
%\put(53,15){\line(3,5){1}}

\end{picture}
\end{figure}

The following ``restricted type'' results will be proved directly.

\begin{theorem}\label{teorema1}
For every vertex $A_i$, $i = 1, \ldots, 6$ there exist
admissible tuples $\alpha$ arbitrarily close to $A_i$ such that the 
form $\Lambda'$ is of restricted
type $\alpha$.
\end{theorem}

\begin{theorem}\label{teorema2}
For the vertices $M_{56}, M_{12}, M_{34}, A$ there exist
admissible tuples $\alpha$ arbitrarily close to them, such that the 
form $\Lambda''$ is of restricted
type $\alpha$.
\end{theorem}

By interpolation of the restricted type estimates 
in the previous Theorem \ref{teorema1}, we first obtain (cf. \cite{mtt:multilinear})

\begin{corollary}\label{restricted1}
Let $\alpha$ be an admissible tuple inside the hexagon $[A_1,...,A_6]$.
Then $\Lambda'$ is of restricted type
$\alpha$.
\end{corollary}
Similarly, we also have

\begin{corollary}\label{restricted2}
Let $\alpha$ be an arbitrary tuple inside the polygon
$[M_{56} M_{34} A M_{12}]$.
Then $\Lambda''$ is of restricted type
$\alpha$.
\end{corollary}

Intersecting these two corollaries we obtain the analogous result for 
$\Lambda$ and the pentagon $[M_{56} M_{34} A_3 A_2 M_{12}]$.

Since one observes that $p_1, p_2, p'_3$ satisfy the hypothesis
of Theorem \ref{main} if and only if 
$(1/p_1,1/p_2,1/p_3)\in [M_{56} M_{34} A_3 A_2 M_{12}]$,
it only remains to convert these restricted type estimates into strong type estimates.  
To do this, one just has to apply (exactly as in \cite{mtt:walshbicarleson}) the multilinear Marcinkiewicz interpolation Theorem 
\cite{janson} in the case of good tuples and the interpolation Lemma 3.11 in \cite{mtt:walshbicarleson} in the case of bad tuples.

This ends the proof of Theorem \ref{main}. Hence, it remains to prove
Theorem \ref{teorema1} and Theorem \ref{teorema2}.

\section{Grids and tiles}\label{grids-sec}

We now start to prove Theorem \ref{teorema1} and Theorem \ref{teorema2}. We look at the symbol $\chi_{\xi_1<\xi_2<N}$ as being the product between
$\chi_{\xi_1<\xi_2}$ and $\chi_{\xi_2<N}$. To construct the multipliers $m'_N$, $m''_N$ and $m'''_N$ we plan to carve these two symbols carefully into smaller pieces.
To take care of the first symbol $\chi_{\xi_1<\xi_2}$, we first need to recall some definitions from \cite{mtt:fourierbiest}.

\begin{definition}\label{mesh-def}
For $n\geq 1$ and $\sigma \in \{ 0, \frac{1}{3}, \frac{2}{3} \}^n$.  We define the \emph{shifted $n$-dyadic mesh} $D = D^n_{\sigma}$ to be the collection of cubes of the form
$$ D^n_{\sigma} := \{ 2^j (k + (0,1)^n + (-1)^j \sigma)| j \in \Z, \quad k \in \Z^n \}.$$
We define a \emph{shifted dyadic cube} to be any member of a shifted $n$-dyadic mesh.
\end{definition}

Observe that for every cube $Q$, there exists a shifted dyadic cube $Q'$ such that $Q \subseteq \frac{7}{10} Q'$ and $|Q'| \sim |Q|$.

\begin{definition}\label{sparse-grid}
A subset $D'$ of a shifted $3$-dyadic grid $D$ is called sparse, if for any two 
cubes $Q,Q'$ in $D$ with $Q\neq Q'$ we have $|Q|<|Q'|$ implies $|10^{9}Q|<|Q'|$ and 
$|Q|=|Q'|$ implies $10^9Q\cap 10^9 Q'=\emptyset$.
\end{definition}

Observe that any subset of a shifted $3$-dyadic grid, can be 
split into $O(1)$ sparse subsets.

\begin{definition}\label{cube-rank-def}  Let $\sigma \in \{ 0, \frac{1}{3}, \frac{2}{3} \}^3$ be a shift.  A collection $\Q \subset D^3_\sigma$ of cubes is said to have \emph{rank 1} if one has the following properties for all $Q, Q' \in \Q$:
\begin{itemize}
\item If $Q \neq Q'$, then $Q_i \neq Q'_i$ for all $i=1,2,3$ (in particular, the $Q$ are disjoint).
\item If $3Q'_j \subset 3Q_j$ for some $j=1,2,3$, then $10^7Q'_i \subset 10^7Q_i$ for all $1 \leq i \leq 3$.  
\item If we further assume that $|Q'| < |10^9Q|$, then we have $3Q'_i \cap 3Q_i = \emptyset$ for all $i \neq j$.
\end{itemize}
\end{definition}
Also, if $Q$ is a cube in $\R^3$, denote by $-Q$ the reflected 
cube about the origin; similarly for intervals in $\R$.

\begin{definition}  Let $\sigma = (\sigma_1,\sigma_2,\sigma_3) \in \{ 0, \frac{1}{3}, \frac{2}{3} \}^3$, and let $1 \leq i \leq 3$.  An $i$-tile with shift $\sigma_i$ is a rectangle $P = I_P \times \omega_P$ with area 1 and with $I_P \in D^1_0$, $\omega_P \in D^1_{\sigma_i}$.  A tri-tile with shift $\sigma$ is a $3$-tuple 
$\pv = (P_1, P_2, P_3)$ such that each $P_i$ is an $i$-tile with shift $\sigma_i$, and the $I_{P_i} = I_\pv$ are independent 
of $i$.  The frequency cube $Q_\pv$ of a tri-tile is defined to be
$Q_\pv= \prod_{i=1}^3 \omega_{P_i}$.
\end{definition}

We shall sometimes refer to $i$-tiles with shift $\sigma$ just as $i$-tiles, or even as tiles, if the parameters $\sigma$, $i$ are unimportant.

\begin{definition} A set $\Pv$ of tri-tiles is called sparse, if all tri-tiles in $\Pv$
have the same shift and the set $\{Q_{\pv}:\pv\in \Pv\}$ is sparse.
\end{definition}

Again, any set of tri-tiles can be split into $O(1)$ sparse subsets.

\begin{definition}  Let $P$ and $P'$ be tiles.  We define a ``relaxed ordering'' and write $P' < ^r P$ if $I_{P'} \subsetneq I_P$ and $3\omega_P \subseteq 3\omega_{P'}$, 
and $P' \leq ^r P$ if $P' < ^r P$ or $P' = P$.
We write $P' \lesssim ^r P$ if $I_{P'} \subseteq I_P$ and $10^7\omega_P \subseteq 10^7 \omega_{P'}$.  We write $P' \lesssim ^{'r} P$ if $P' \lesssim ^r P$ and
$P' \not \leq ^r P$.
\end{definition}

This ``relaxed ordering'' $< ^r$ is in the spirit of that in
Fefferman \cite{fefferman} or Lacey and Thiele \cite{laceyt},  but slightly different as $P'$ and $P$ do not quite have to intersect.  
This is more convenient for technical purposes.

\begin{definition}\label{rank-def}  A collection $\Pv$ of tri-tiles is said to have \emph{rank 1} if one has the following properties for all $\pv, \pv' \in \Pv$:
\begin{itemize}
\item If $\pv \neq \pv'$, then $P_j \neq P'_j$ for all $j=1,2,3$.
\item If $P'_j \leq ^rP_j$ for some $j=1,2,3$, then $P'_i \lesssim ^rP_i$ for all $1 \leq i \leq 3$.  
\item If we further assume that $|I_{\pv'}| < 10^9 |I_{\pv}|$, then we have $P'_i \lesssim ^{'r} P_i$ for all $i \neq j$.
\end{itemize}
\end{definition}

\begin{definition}\label{packet-def}  Let $P$ be a tile.  A \emph{wave packet on $P$} is a function $\phi_P$ which has Fourier support in $\frac{9}{10} \omega_P$ and obeys the estimates 
\be{decay}
|\phi_P(x)| \lesssim |I_P|^{-1/2} \tilde \chi_I(x)^M
\end{equation}
for all $M > 0$, with the implicit constant depending on $M$.
\end{definition}

Heuristically, $\phi_P$ is $L^2$-normalized and is supported in $P$.

To take care of the second symbol $\chi_{\xi_2<N}$ we need two more definitions.

\begin{definition}\label{grid-def}
Let $u\in [-1,1]$ and $v\in\R$. We define the dyadic grid $D=D_{u,v}$ associated to the parameters $u, v$ to be the collection of all intervals
of the form
$$2^{u}I + v,$$
where $I\in D^1_0$ is an arbitrary dyadic interval on the real line.
\end{definition}

\begin{definition}\label{bitiles-def}
Let $u\in [-1,1]$ and $v\in\R$, as before. A bi-tile associated to the parameters $u, v$ is a rectangle $P=I_P\times\omega_P$ of area $2$ such that
$\omega_P\in D_{u,v}$ and $I_P\in D_{-u, 0}$. To any bi-tile $P$ we define its (sub)-tiles $P_1, P_2$ as being given by
$P_1=I_{P_1}\times\omega_{P_1}$, $P_2=I_{P_2}\times\omega_{P_2}$, where $I_{P_1}=I_{P_2}:=I_P$, $\omega_{P_1}:=\omega_P^l$, $\omega_{P_2}:=\omega_P^r$,
while $\omega_P^l$ and $\omega_P^r$ are the left and right halves of $\omega_P$.
\end{definition}

\section{The symbol $m'_N$}\label{m'-sec}

Fix $N$. In this section we construct the symbol $m'_N$. As in \cite{mtt:fourierbiest}, by a standard partition of unity we can write
$$
\chi_{\xi_1 < \xi_2}(\xi_1,\xi_2,\xi_3) = 
\sum_{\sigma \in \{0, \frac 1 3, \frac 2 3\}^3} \sum_{Q \in \Q_\sigma} \phi_{Q,\sigma}(\xi_1, \xi_2, \xi_3)$$
whenever $\xi_1 + \xi_2 + \xi_3 = 0$, where $\Q_\sigma \subset D^3_\sigma$ is a collection of cubes which intersect the plane $\{\xi_1 + \xi_2 + \xi_3 = 0\}$ and which satisfy the Whitney property
$$ 10^3 \diam(Q) \leq \dist(Q, \{ \xi_1 = \xi_2, \xi_1 + \xi_2 + \xi_3 = 0 \}) \leq 10^5 \diam(Q) $$
for all $Q \in \Q_\sigma$, and for each cube $Q \in \Q_\sigma$, $\phi_{Q,\sigma}$ is a bump function adapted to $\frac{8}{10} Q$.  Note that by refining $\Q_\sigma$ by a finite factor if necessary one can make $\Q_\sigma$ have rank 1 (this refining of $\Q_\sigma$ corresponds to decomposing $\chi_{\xi_1<\xi_2}$ into O(1) pieces).

Also as in \cite{mtt:fourierbiest}, by splitting $\phi_{Q,\sigma}$ as a Fourier series in the $\xi_i$ we can then write
$$
\chi_{\xi_1 < \xi_2}(\xi_1,\xi_2,\xi_3) = \sum_{l \in \Z^3} c_l \sum_{\sigma \in \{0, \frac 1 3, \frac 2 3\}^3} \sum_{Q \in \Q_\sigma} \eta_{Q_1,\sigma,l,1}(\xi_1)
\eta_{Q_2,\sigma,l,2}(\xi_2)
\eta_{Q_3,\sigma,l,3}(\xi_3),$$
where $c_l$ is a rapidly decreasing sequence 
and $\eta_{Q_j, \alpha, l, j}$ is a bump function adapted to $\frac{9}{10} Q_j$ uniformly in $l$.

Let now $\phi$ be a Schwartz function so that $\widehat{\phi}$ is supported on $[-1/8, 1/8]$ and which equals $1$ on
$[-1/16, 1/16]$. If $\omega$ is a dyadic interval, let $\phi_{\omega}$ be the function defined via the Fourier transform by

\begin{equation}
\widehat{\phi_{\omega}}:= T_{c(\omega)} D^{\infty}_{|\omega|}\widehat{\phi}
\end{equation}
where in general
$$T_y f(x):= f(x-y)$$
$$D^{\infty}_{\lambda} f(x):= f(\lambda^{-1}x)$$
are the translation and dilation operators respectively. As in \cite{laceyt1} (see also \cite{fefferman}), one can write $\chi_{\xi <N}$ as

\begin{equation}\label{simbol2}
\chi_{\xi <N}= c\cdot\lim_{n\rightarrow \infty} \frac{1}{|R_n|}
\int_{R_n}\sum_{\omega\in D_{0,0}}\chi_{\omega^r}(2^{-k}(N+\eta)) T_{-\eta} D^{\infty}_{2^k} \widehat{\phi_{\omega^l}}(\xi) d\eta dk
\end{equation}

$$:=\int{\!\!\!\!\!\! -} \sum_{\omega\in D_{0,0}}\chi_{\omega^r}(2^{-k}(N+\eta)) T_{-\eta} D^{\infty}_{2^k} \widehat{\phi_{\omega^l}}(\xi) d\eta dk,$$
where $R_n:= [0,1]\times [-n, n]\subseteq \R^2$ and $\omega^l, \omega^r$ are the left and right halvess of $\omega$. In particular, this implies that
on the line $\xi+\xi'=0$ in $\R^2$ one has

$$\chi_{\xi <N}(\xi,\xi')=
\int{\!\!\!\!\!\! -} \sum_{\omega\in D_{0,0}}\chi_{\omega^r}(2^{-k}(N+\eta)) T_{-\eta} D^{\infty}_{2^k} \widehat{\phi_{\omega^l}}(\xi) 
T_{-\eta} D^{\infty}_{2^k} \widehat{\tilde{\phi}_{-\omega^l}}(\xi')
d\eta dk
$$

$$=\int{\!\!\!\!\!\! -}\sum_{\omega\in D_{k,\eta}}\chi_{\omega^r}(N)\widehat{\phi_{\omega^l}}(\xi) \widehat{\tilde{\phi}_{-\omega^l}}(\xi') d\eta dk,$$
where this time $\tilde{\phi}$ is a Schwartz function, such that $\widehat{\tilde{\phi}}$ is supported on $[-1/4, 1/4]$ and equals
$1$ on $[-1/8, 1/8]$.
As a consequence, the expression

\begin{equation}\label{E}
E_N(\xi_1, \xi_2, \xi_3):= \sum_{l \in \Z^3} c_l \sum_{\sigma \in \{0, \frac 1 3, \frac 2 3\}^3} \int{\!\!\!\!\!\! -}\sum_{Q \in \Q_\sigma} 
\sum_{\omega\in D_{k,\eta}; -Q_3\subseteq\omega^l}
\eta_{Q_1,\sigma,l,1}(\xi_1)
\eta_{Q_2,\sigma,l,2}(\xi_2)
\eta_{Q_3,\sigma,l,3}(\xi_3)\cdot
\end{equation}
$$\cdot\chi_{\omega^r}(2N)\widehat{\phi_{\omega^l}}(\xi_1+\xi_2) \widehat{\tilde{\phi}_{-\omega^l}}(\xi_3) d\eta dk$$
on the hyperplane $\xi_1+\xi_2+\xi_3=0$ is equal to $\chi_{\xi_1<\xi_2} \chi_{\frac{\xi_1+\xi_2}{2}<N}$ when
$|\xi_1-\xi_2|\ll |\frac{\xi_1+\xi_2}{2}-N|$ (under the latter constraint the condition $-Q_3\subseteq \omega^l$ is automatic for nonzero summands).
Then, we simply define the symbol $m'_N$ by setting $m'_N:= E_N(\xi_1, \xi_2, -\xi_1-\xi_2)$ (in other words, as we said before, the intuitive description
of $m'_N$ is a smooth restriction of $\chi_{\xi_1<\xi_2<N}$ to the region where $\xi_1$ and $\xi_2$ are closer to each other than they are to $N$).

\section{Discretization of $T'$}

Recall that the bilinear operator $T'$ has been defined by the formula

$$T'(f_1, f_2)(x):= T_{m'_{N(x)}}(f_1, f_2)(x).$$
The aim of the present section is to reduce Theorem \ref{teorema1} to the following discretized version of it.

\begin{theorem}\label{discrete}
Let $\sigma\in \{ 0, \frac{1}{3}, \frac{2}{3} \}^3$ be a shift and let $u\in [0,1]$, $v\in\R$.
Let $\P$ be a finite collection of bi-tiles associated to the parameters $u,v$ and $\vec{\Q}$ be a finite collection of tri-tiles of rank $1$ with shift $\sigma$.
For each $i=1,2$ and $P \in \P$, let $\phi_{P_i} = \phi_{P_i,i}$ and  $\tilde\phi_{P_i} = \tilde\phi_{P_i,i}$  be wave packets on $P_i$.  Similarly for each $i=1,2,3$ 
and $\qv \in \Qv$ let $\tilde{\tilde\phi}_{Q_i} = \tilde{\tilde\phi}_{Q_i,i}$ be a wave packet on $Q_i$.  
Define the form $\Lambda'_{\P,\Qv}$ by
$$\Lambda'_{\P,\Qv}(f_1,f_2,f_3) :=
\sum_{P \in \P} 
\langle B_{P_1,\vec{\Q}}(f_1, f_2), \phi_{P_1} \rangle
\langle \tilde\phi_{P_1}\chi_{\{x : N(x)\in\omega_{P_2}\}}, f_3 \rangle
$$
where
$$
B_{P_1,\vec{\Q}}(f_1, f_2)  := \sum_{\qv \in \Qv: \omega_{Q_3} \subseteq \omega_{P_1}}
\frac{1}{|I_\qv|^{1/2}}
\langle f_1,  \tilde{\tilde\phi}_{Q_1} \rangle
\langle f_2,  \tilde{\tilde\phi}_{Q_2} \rangle
\tilde{\tilde\phi}_{Q_3}.$$
Then for every vertex $A_j$ $j=1,...,6$, there exist admissible tuples $\alpha$ arbitrarily close to $A_j$ such that the form $\Lambda'_{\P,\Qv}$ is of restricted type
$\alpha$, uniformly in the parameters $\sigma$, $u, v$ $\P$, $\Qv$, $\phi_{P_i}$, $\tilde\phi_{P_i}$      $\tilde{\tilde\phi}_{Q_i}$. 
Furthermore, in the case that $\alpha$ has a bad index $j$, the restricted type is uniform in the sense that the  major subset $E'_j$ can be chosen independently of 
the parameters just mentioned.
\end{theorem}

\setlength{\unitlength}{0.8mm}
\begin{figure}\label{tiles'}
\caption{The tiles of $\Lambda'$}
\begin{picture}(190,100)

\put(150,5){\line(0,1){80}}
\put(155,5){\line(0,1){80}}
\put(150,5){\line(1,0){5}}
\put(150,45){\line(1,0){5}}
\put(150,85){\line(1,0){5}}

\put(40,10){\line(1,0){80}}
\put(40,12.5){\line(1,0){80}}
\put(40,20){\line(1,0){80}}
\put(40,22.5){\line(1,0){80}}
\put(40,30){\line(1,0){80}}
\put(40,32.5){\line(1,0){80}}

\put(40,10){\line(0,1){2.5}}
\put(40,20){\line(0,1){2.5}}
\put(40,30){\line(0,1){2.5}}

\put(120,10){\line(0,1){2.5}}
\put(120,20){\line(0,1){2.5}}
\put(120,30){\line(0,1){2.5}}

\put(157,25){$P_1$}
\put(157,65){$P_2$}

\put(31,10){$Q_1$}
\put(31,20){$Q_2$}
\put(31,30){$Q_3$}

\end{picture}
\end{figure}
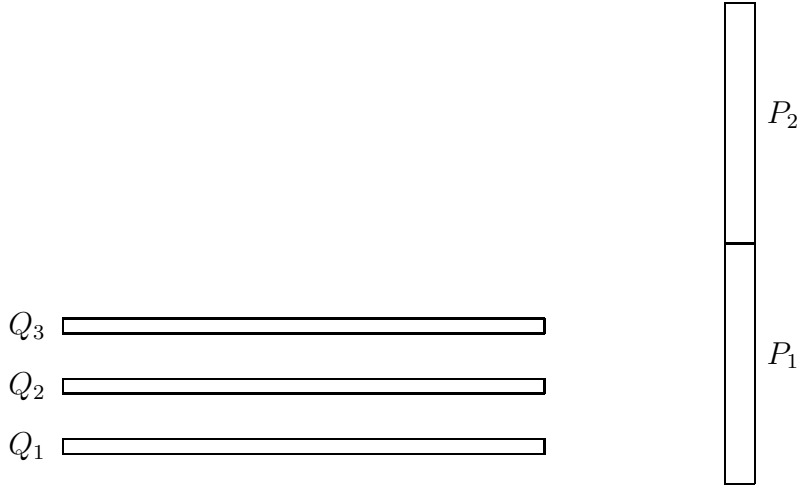

In the rest of this section we show how Theorem \ref{teorema1} can be deduced from Theorem \ref{discrete}.  
To see this, we  just have to calculate the form $\Lambda'$ carefully. We fix $N$ and first look at the trilinear form associated to the operator $T_{m'_N}$.
It is given by the formula

$$\int_{\R}T_{m'_N}(f_1, f_2)(x) f_3(x) dx = \int\delta(\xi_1+\xi_2+\xi_3)E_N(\xi_1, \xi_2, \xi_3)
\widehat{f}_1(\xi_1)\widehat{f}_2(\xi_2)\widehat{f}_3(\xi_3)
d\xi_1d\xi_2d\xi_3,$$
where $\delta$ is the Dirac delta function. Clearly, $E_N$ is an average of some simpler expressions depending on the parameters $l, \sigma, \eta, k$. 
We fix all of them and look at the corresponding form given now by

\begin{equation}\label{simpla}
\int\delta(\xi_1+\xi_2+\xi_3)
\sum_{Q,\omega; -Q_3\subseteq\omega^l}
\eta_{Q_1}(\xi_1)
\eta_{Q_2}(\xi_2)
\eta_{Q_3}(\xi_3)
\chi_{\omega^r}(2N)\widehat{\phi_{\omega^l}}(\xi_1+\xi_2) \widehat{\tilde{\phi}_{-\omega^l}}(\xi_3)\cdot
\end{equation}

$$\cdot\widehat{f}_1(\xi_1)\widehat{f}_2(\xi_2)\widehat{f}_3(\xi_3)
d\xi_1d\xi_2d\xi_3.$$

To calculate further (\ref{simpla}), let us first consider the simpler object
$$
\left<f_3,a_Q\right>
:=\int \delta(\xi_1+\xi_2+\xi_3) 
\eta_{Q_1}(\xi_1) \eta_{Q_2}(\xi_2) \eta_{Q_3}(\xi_3)
\widehat{f}_1(\xi_1)\widehat{f}_2(\xi_2)\widehat{f}_3(\xi_3)
d\xi_1d\xi_2d\xi_3$$
for some cube $Q \in \Q$.  By Plancherel this is equal to
\bas
&
\int (f_1 * \check \eta_{Q_1})(x)
(f_2 * \check \eta_{Q_2})(x)
(f_3 * \check \eta_{Q_3})(x)
\ dx\\
&= l(Q)^{3/2}
\int \langle f_1, \phi_{Q_1,x} \rangle
\langle f_2, \phi_{Q_2,x} \rangle
\langle f_3, \phi_{Q_3,x} \rangle\ dx
\end{align*}
where
$$ \phi_{Q_j,x}(y) := l(Q)^{-1/2} \overline{\check \eta_{Q_j}(x-y)}$$
and $l(Q)$ is the side-length of $Q$.  We can rewrite this as
$$
\int_0^1  \sum_{\pv: Q_\pv = Q} |I_\pv|^{-1/2}
\langle f_1, \phi_{P_1,t,1} \rangle
\langle f_2, \phi_{P_2,t,2} \rangle
\langle f_3, \phi_{P_3,t,3} \rangle
\ dt
$$
where $\pv$ ranges over all tri-tiles with frequency cube $Q$ and spatial interval $I_\pv$ in $D^1_0$, $\phi_{P_j,t,j}$ is the function
$$ \phi_{P_j,t,j} := \phi_{Q_j,x_{\pv} + |I_\pv| t}$$
and $x_{\pv}$ is the center of $I_\pv$.  Note that $\phi_{P_j,t,j}$ is a wave packet on $P_j$ uniformly in $t$.

Similarly, we consider

$$\left<f, b_{\omega}\right>:=
\int \delta(\xi+\xi')\chi_{\omega^r}(2N) \widehat{\phi}_{\omega^l}(\xi)\widehat{\tilde\phi}_{-\omega^l}(\xi') d\xi d\xi'$$
for some fixed $\omega$ as above. Again, by Plancherel, this is equal to

$$\int_{\R}(f*\phi_{\omega^l})(x) (f_3\chi_{\omega^r}(2N)*\tilde\phi_{-\omega^l})(x) dx$$

$$=\int_0^1\sum_{P; \omega_P=\omega} 
\langle f, \phi_{P_1,t,1} \rangle
\langle f_3\chi_{\omega_{P_2}}(2N)  , \tilde\phi_{\overline{P}_1,t,1} \rangle dt$$
where $P$ ranges over all bi-tiles with frequency interval $\omega$, $\overline{P}_1:= I_{P_1}\times (-\omega_{P_1})$ and

$$\phi_{P_1,t,1}(y):=|I_P|^{1/2}\overline{\phi_{\omega^l}(x_P + |I_P|t - y)},$$

$$\phi_{\overline{P}_1,t,1}(y):=|I_P|^{1/2}\overline{\tilde\phi_{-\omega^l}(x_P + |I_P|t - y)}$$
where $x_P$ is the center of the spatial interval $I_P$.
As a consequence, our expression in (\ref{simpla}) can be written as
$$\sum_{Q,\omega; -Q_3\subseteq\omega^l}
\int \overline{\widehat{a_Q}(\tau)\widehat{b_{\omega}}(-\tau)} \, d\tau$$

$$=\int_0^1\int_0^1\sum_{P\in \P}
\langle B_{P_1, t'}(f_1, f_2), \phi_{P_1,t,1} \rangle
\langle f_3\chi_{\omega_{P_2}}(2N)  , \tilde\phi_{\overline{P}_1,t,1} \rangle dt dt'$$
where

$$B_{P_1,t'}(f_2,f_3) := \sum_{\qv: Q_\qv \in \Q, - \omega_{Q_3} \subseteq \omega_{P_1}}
\frac{1}{|I_\qv|^{1/2}}
\langle f_1, \phi_{Q_1,t',1} \rangle
\langle f_2, \phi_{Q_2,t',2} \rangle
\overline{\phi_{Q_3,t',3}}.$$
Note that the collection $\Qv$ of tri-tiles $\qv$ has rank 1.  Observe that we can get rid of the complex conjugation sign in the definition of $B_{P_1,t'}$ by redefining 
$Q_3$ to be $-Q_3$ and 
redefining $\phi_{Q_3,t',3}$ accordingly; 
this 
also replaces the condition
$- \omega_{Q_3} \subseteq \omega_{P_1}$ by the condition $\omega_{Q_3} \subseteq \omega_{P_1}$, but it does not change the rank one property of the collection $\Qv$.

This means that the operator corresponding to the form (\ref{simpla}) becomes now

$$\int_0^1\int_0^1
\langle B_{P_1, t'}(f_1, f_2), \phi_{P_1,t,1} \rangle
\overline{\tilde\phi}_{\overline{P}_1,t,1} \chi_{\omega_{P_2}}(2N) dt dt'$$

$$=\int_0^1\int_0^1
\langle B_{P_1, t'}(f_1, f_2), \phi_{P_1,t,1} \rangle
\tilde\phi_{P_1,t,1} \chi_{\omega_{P_2}}(2N) dt dt'.$$
In particular, our operator $T'$ is in fact an average of simpler operators $U$ of the form

$$U(f_1, f_2)(x) = \int_0^1\int_0^1
\langle B_{P_1, t'}(f_1, f_2), \phi_{P_1,t,1} \rangle
\tilde\phi_{P_1,t,1}(x) \chi_{\omega_{P_2}}(2N(x)) dt dt'.$$

The claim then follows by integrating the conclusion of Theorem \ref{discrete} over $t$, $t'$, $\eta$, $k$, summing over $\sigma$, $l$ and
 using the uniformity assumptions of that Theorem.  
(The finiteness condition on $\P$ and $\Qv$ can be removed by the usual limiting arguments.)

\section{The symbol $m''_N$}

Fix $N$. In this section we construct the symbol $m''_N$. First, by using (\ref{simbol2}) we know that

$$\chi_{\xi_1 < N} \cdot \chi_{\xi_2 < N}=$$

$$=\int{\!\!\!\!\!\! -}\int{\!\!\!\!\!\! -}\sum_{\omega, \omega'\in D_{0,0}}\chi_{\omega^r}(2^{-k}(N+\eta)) T_{-\eta} D^{\infty}_{2^k} \widehat{\phi_{\omega^l}}(\xi_1) 
\chi_{\omega^{'r}}(2^{-k'}(N+\eta')) T_{-\eta'} D^{\infty}_{2^{k'}} \widehat{\phi_{\omega^{'l}}}(\xi_2) d\eta d\eta' dk dk'=$$

$$\int{\!\!\!\!\!\! -}\int{\!\!\!\!\!\! -}\sum_{\omega\in D_{k,\eta}}\sum_{\omega'\in D_{k',\eta'}}
\chi_{\omega^r}(N)\widehat{\phi_{\omega^l}}(\xi_1)
\chi_{\omega^{'r}}(N)\widehat{\phi_{\omega^{'l}}}(\xi_2) d\eta d\eta' dk dk'.$$

Then, we define $m''_N$ to be given by the following formula

$$m''_N:= \int{\!\!\!\!\!\! -}\int{\!\!\!\!\!\! -}\sum_{\omega^{'r}\cap\omega^r\neq\emptyset ; |\omega^{'r}|\leq |\omega^r|}
\chi_{\omega^r}(N)\widehat{\phi_{\omega^l}}(\xi_1)
\chi_{\omega^{'r}}(N)\widehat{\phi_{\omega^{'l}}}(\xi_2) d\eta d\eta' dk dk'.$$
We observe that the above expression is equal to $\chi_{\xi_1 < \xi_2 < N}$ on the cone $|\xi_2-N| \ll |\xi_1-N|$ (under the latter constraint,
the condition $\omega^{'r}\cap\omega^r\neq\emptyset ; |\omega^{'r}|\leq |\omega^r|$ is automatic for nonzero summands). In the end we set
$m'''_N:= \chi_{\xi_1 < \xi_2 < N} - m'_N - m''_N$ and remark that all the requirements described in Section 1 are satisfied.

\section{Discretization of $T''$}

Recall that the bilinear operator $T''$ has been defined by the formula

$$T''(f_1, f_2)(x):= T_{m''_{N(x)}}(f_1, f_2)(x).$$
The aim of the present section is to reduce Theorem \ref{teorema2} to the following discretized version of it.

\begin{theorem}\label{discrete1}
Let  $u, u'\in [0,1]$, $v, v'\in\R$.
Let $\P$, $\Q$ be finite collections of bi-tiles associated to the parameters $u, v$ and $u', v'$ respectively. 
For each $i=1,2$ and $P \in \P$, let $\phi_{P_i} = \phi_{P_i,i}$ and  $\tilde\phi_{P_i} = \tilde\phi_{P_i,i}$  be wave packets on $P_i$.  Similarly for each $i=1,2$ 
and $Q\in \Q$ let $\tilde{\tilde\phi}_{Q_i} = \tilde{\tilde\phi}_{Q_i,i}$ and $\tilde{\tilde{\tilde\phi}}_{Q_i} = \tilde{\tilde{\tilde\phi}}_{Q_i,i}$ be  wave packets 
on $Q_i$.  
Define the form $\Lambda''_{\P,\Q}$ by
$$\Lambda''_{\P,\Q}(f_1,f_2,f_3) :=
\sum_{P \in \P} 
\langle f_1, \phi_{P_1} \rangle
\langle \tilde\phi_{P_1} \chi_{\{x : N(x)\in \omega_{P_2}\}} C_{P_2, \Q}(f_2), f_3\rangle
$$
where
$$C_{P_2, \Q}(f_2):= \sum_{\omega_{Q_2}\cap\omega_{P_2}\neq \emptyset ; |\omega_{Q_2} | < |\omega_{P_2} | }
\langle f_2, \tilde{\tilde\phi}_{Q_1} \rangle
\tilde{\tilde{\tilde\phi}}_{Q_1} \chi_{\{x : N(x)\in \omega_{Q_2}\}}
.$$
Then for the vertices $M_{56}$, $M_{12}$, $A_2$, there exist admissible tuples $\alpha$ arbitrarily close to them such that the form $\Lambda''_{\P,\Q}$ is of restricted type
$\alpha$, uniformly in the parameters $\sigma$, $u, v$ $\P$, $\Q$, $\phi_{P_i}$, $\tilde\phi_{P_i}$,  $\tilde{\tilde\phi}_{Q_i}$, $\tilde{\tilde{\tilde\phi}}_{Q_i}$ 
Furthermore, in the case that $\alpha$ has a bad index $j$, the restricted type is uniform in the sense that the  major subset $E'_j$ can be chosen independently of 
the parameters just mentioned.
\end{theorem}

\setlength{\unitlength}{0.8mm}
\begin{figure}\label{tiles''}
\caption{The tiles of $\Lambda''$}
\begin{picture}(190,100)

\put(150,5){\line(0,1){80}}
\put(155,5){\line(0,1){80}}
\put(150,5){\line(1,0){5}}
\put(150,45){\line(1,0){5}}
\put(150,85){\line(1,0){5}}

\put(40,60){\line(1,0){80}}
\put(40,62.5){\line(1,0){80}}
\put(40,65){\line(1,0){80}}
\put(40,60){\line(0,1){5}}
\put(120,60){\line(0,1){5}}

\put(157,25){$P_1$}
\put(157,65){$P_2$}

\put(31,55){$Q_1$}
\put(31,70){$Q_2$}

\end{picture}
\end{figure}
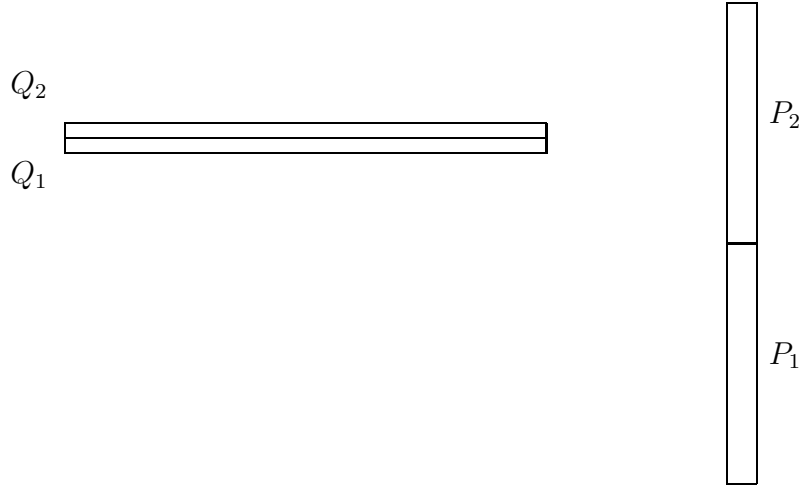

In the rest of this section we show how Theorem \ref{teorema2} can be deduced from Theorem \ref{discrete1}.  
To see this, as before, we just have to calculate the form $\Lambda''$ carefully. We fix $N$ and first look at the trilinear form associated to the operator $T_{m''_N}$.
Clearly, $m''_N$ is an average of some simpler symbols depending on the parameters $\eta$, $\eta'$, $k$, $k'$. We fix all of 
these parameters and look at the corresponding operator given by

\begin{equation}
\int\sum_{\omega^{'r}\cap\omega^r\neq\emptyset ; |\omega^{'r}|\leq |\omega^r|}
\chi_{\omega^r}(N)\widehat{\phi_{\omega^l}}(\xi_1)
\chi_{\omega^{'r}}(N)\widehat{\phi_{\omega^{'l}}}(\xi_2)
\widehat{f_1}(\xi_1)\widehat{f_2}(\xi_2) e^{2\pi i x(\xi_1 + \xi_2)} d\xi_1 d\xi_2
\end{equation}

$$=\sum_{\omega^{'r}\cap\omega^r\neq\emptyset ; |\omega^{'r}|\leq |\omega^r|}
(f_1*\phi_{\omega^l})(x)\chi_{\omega^r}(N)
(f_2*\phi_{\omega^{'l}})(x)\chi_{\omega^{'r}}(N)$$
which can be discretized as before into

$$\int_0^1\int_0^1
\sum_{\omega_{Q_2}\cap\omega_{P_2}\neq \emptyset ; |\omega_{Q_2} | < |\omega_{P_2} |  }
\langle f_1, \phi_{P_1, t, 1} \rangle \tilde\phi_{P_1, t, 1} \chi_{\omega_{P_2}}(N)
\langle f_2, \phi_{Q_1, t', 1} \rangle \tilde\phi_{Q_1, t', 1} \chi_{\omega_{Q_2}}(N)
dt dt'.$$
The claim then follows using a similar argument.

\section{trees}

The standard approach to prove the desired estimates for the forms $\Lambda'$ and $\Lambda''$
is to organize our fixed collections of tri-tiles and bi-tiles into trees as in \cite{fefferman}.
We may assume and shall do so for the rest of this article that our collections of tri-tiles are sparse.

Firstly, we define trees in the context of of tri-tiles, by using the ``relaxed ordering'' $\leq ^r$ considered before.

\begin{definition} Let $\vec{\P}$ be a collection of tri-tiles. For any $1 \leq j \leq 3$ and a tri-tile $\pv_T \in \Pv$, define a $j$-tree with top $\pv_T$ to be a collection of tri-tiles $T \subseteq \Pv$ such that
$$ P_j \leq^r P_{T,j} \hbox{ for all } \pv \in T,$$
where $P_{T,j}$ is the $j$ component of $\pv_T$.  We write $I_T$ and $\omega_{T,j}$ for $I_{\pv_T}$ and $\omega_{P_{T,j}}$ respectively.  We say that $T$ is a tree if it is a $j$-tree for some $1 \leq j \leq 3$.
\end{definition}

Note that $T$ does not necessarily have to contain its top $\pv_T$.

\begin{definition} Let $1 \leq i \leq 3$.  Two trees $T$, $T'$ are said to be \emph{strongly $i$-disjoint} if 
\bi
\item $P_i \neq P'_i$ for all $\pv \in T$, $\pv' \in T'$.
\item Whenever $\pv \in T$, $\pv' \in T'$ are such that
$2\omega_{P_i} \cap 2\omega_{P'_i}\neq \emptyset$, then one has
$I_{\pv'} \cap I_T = \emptyset$, and similarly with $T$ and $T'$ reversed.
\end{itemize}
\end{definition}

Note that if $T$ and $T'$ are strongly $i$-disjoint, 
then $I_P\times 2\omega_{P_i} \cap I_{P'}\times 2\omega_{{P'}_i} = \emptyset$ 
for all $\pv \in T$, $\pv' \in T'$.

Given that $\vec{\P}$ is sparse, it is easy to see that
if $T$ is an $i$- tree, then for all $\pv,\pv'\in T$ and $j\neq i$ we have
$$ \omega_{P_j} = \omega_{P'_j}$$ or 
$$ 2\omega_{P_j} \cap 2\omega_{P'_j}=\emptyset.$$

Secondly, we define trees in the context of bi-tiles, this time by using the ``classical ordering'' $\leq ^c$ as in \cite{fefferman}.
\begin{definition}
We define a partial ordering on the set of bi-tiles $\P$ by $P\leq ^c P'$ if $I_P\subseteq I_{P'}$ and $\omega_{P'}\subseteq \omega_P$.
A set $T$ of bi-tiles is called a tree, if there is a tile $P_T$, the top of the tree, such that $P\leq ^c P_T$ for all $P\in T$.
If $j=1,2$, a tree is called $j$-tree if $\omega_{P_{T_j}}\subseteq \omega_{P_j}$ for all $P\in T$. Notice that as before, we do not require the top 
to be an element of the tree.
\end{definition}

\section{Tile norms}
We start the study of the form $\Lambda'_{\P, \Qv}$.
In this section we mostly recall some definitions and results from \cite{mtt:fourierbiest} and \cite{mtt:walshbicarleson}.
In the paper we shall be frequently estimating expressions of the form
\be{trilinear}
| \sum_{\pv \in \Pv} \frac{1}{|I_\pv|^{1/2}} a^{(1)}_{P_1} a^{(2)}_{P_2} a^{(3)}_{P_3}|
\end{equation}
where $\Pv$ is a collection of tri-tiles and $a^{(j)}_{P_j}$ are complex numbers for $\pv \in \Pv$ and $j = 1,2,3$.  In some cases (e.g. if one only wished to treat the Bilinear Hilbert transform) we just have
\be{aj-def}
a^{(j)}_{P_j} = \langle f_j, \phi_{P_j} \rangle
\end{equation}
 but we will have more sophisticated sequences $a^{(j)}_{P_j}$ when dealing with $\Lambda'_{\P,\Qv}$.  

In \cite{mtt:walshbiest} the following (standard) norms on sequences of tiles were introduced: 

\begin{definition}\label{size-def}
Let $\Pv$ be a finite collection of tri-tiles, $j=1,2,3$, and let $(a_{P_j})_{\pv \in \Pv}$ be a sequence of complex numbers.  We define the \emph{size} of this sequence by
$$ \size_j( (a_{P_j})_{\pv \in \Pv} ) := \sup_{T \subset \Pv}
(\frac{1}{|I_T|} \sum_{\pv \in T} |a_{P_j}|^2)^{1/2}$$
where $T$ ranges over all trees in $\Pv$ which are $i$-trees for some $i \neq j$.
We also define the \emph{energy} of the sequence by
$$ \energy_j((a_{P_j})_{\pv \in \Pv} ) := \sup_{\D \subset \Pv}
(\sum_{\pv \in \D} |a_{P_j}|^2)^{1/2}$$
where $\D$ ranges over all subsets of $\Pv$ such that the tiles $\{ P_j: \pv \in \D \}$ are pairwise disjoint.
\end{definition}

The size measures the extent to which the sequence $a_{P_j}$ can concentrate on a single tree and should be thought of as a phase-space variant of the BMO norm.  The energy is a phase-space variant of the $L^2$ norm.  As the  notation suggests, the number $a_{P_j}$ should be thought of as being associated with the tile $P_j$ rather than the full tri-tile $\pv$.

In the Walsh model the energy is a tractable quantity; for instance, if $a_{P_j}$ is given by \eqref{aj-def} then one can control the energy by $\|f_j\|_2$ thanks to the perfect orthogonality of the Walsh wave packets.  However, in the Fourier case the orthogonality is too poor to give a usable bound on the energy, and so we must instead use a more technical substitute.

\begin{definition}\label{energy-def}
Let the notation be as in Definition \ref{size-def}.  We define the \emph{modified energy} of the sequence $(a_{P_j})_{\pv \in \Pv}$ by
\be{energy-mod-def}
\modenergy_j((a_{P_j})_{\pv \in \Pv} ) := \sup_{n \in \Z}
\sup_{\T}
2^{n} (\sum_{T \in \T} |I_T|)^{1/2}
\end{equation}
where $\T$ ranges over all collections of strongly $j$-disjoint trees in $\Pv$ such that
$$ (\sum_{\pv \in T} |a_{P_j}|^2)^{1/2} \ge 2^n |I_T|^{1/2}$$
for all $T \in \T$, and
$$ (\sum_{\pv \in {T'}} |a_{P_j}|^2)^{1/2} \le 2^{n+1} |I_{T'}|^{1/2}$$
for all sub-trees $T' \subset T \in \T$.
\end{definition}

The reader may easily verify that the modified energy is always dominated by the energy, and that we have the monotonicity property
$$
\modenergy_j((a_{P_j})_{\pv \in \Pv'} ) \leq
\energy_j((a_{P_j})_{\pv \in \Pv} )$$
whenever $\Pv' \subset \Pv$. 

The usual BMO norm can be written using an $L^2$ oscillation or an $L^1$ oscillation, and the two notions are equivalent thanks to the John-Nirenberg inequality.  The analogous statement for size is

\begin{lemma}\label{jn}
Let $\Pv$ be a finite collection of tri-tiles, $j=1,2,3$, and let $(a_{P_j})_{\pv \in \Pv}$ be a sequence of complex numbers. Then
\be{jn-est}
\size_j( (a_{P_j})_{\pv \in \Pv} ) \sim \sup_{T \subset \Pv}
\frac{1}{|I_T|} \| ( \sum_{\pv \in T} |a_{P_j}|^2 \frac{\chi_{I_\pv}}{|I_\pv|} )^{1/2} \|_{L^{1,\infty}(I_T)}
\end{equation}
where $T$ ranges over all trees in $\Pv$ which are $i$-trees for some $i \neq j$.
Similarly, for $1<p<\infty$, one has

$$\size_j( (a_{P_j})_{\pv \in \Pv} ) \sim \sup_{T \subset \Pv}
\frac{1}{|I_T|^{1/p}} \| ( \sum_{\pv \in T} |a_{P_j}|^2 \frac{\chi_{I_\pv}}{|I_\pv|} )^{1/2} \|_{L^{p}(I_T)},$$
and the implicit constants are allowed to depend on $p$.
\end{lemma}

\begin{proof}
The  same as in \cite{mtt:walshbiest}, Lemma 4.2 .
\end{proof}

The following estimate is standard, see \cite{mtt:fourierbiest}, Proposition 6.5.  
This is the main combinatorial tool needed to obtain estimates on \eqref{trilinear}.

\begin{proposition}\label{abstract}
Let $\Pv$ be a finite collection of tri-tiles, and for each $\pv \in \Pv$ and $j=1,2,3$ let $a^{(j)}_{P_j}$ be a complex number.  Then
\be{trilinear-est}
| \sum_{\pv \in \Pv} \frac{1}{|I_\pv|^{1/2}} a^{(1)}_{P_1} a^{(2)}_{P_2} a^{(3)}_{P_3}|
\lesssim \prod_{j=1}^3 
\size_j( (a^{(j)}_{P_j})_{\pv \in \Pv} )^{\theta_j}
\modenergy_j( (a^{(j)}_{P_j})_{\pv \in \Pv} )^{1-\theta_j}
\end{equation}
for any $0 \leq \theta_1, \theta_2, \theta_3 < 1$ with $\theta_1 + \theta_2 + \theta_3 = 1$, with the implicit constant depending on the $\theta_j$.
\end{proposition}

Note that this Proposition is stronger than that of the corresponding statement (\cite{mtt:walshbiest}, Proposition 4.3) for the unmodified energy.  
Of course, in order to use Proposition \ref{abstract} we will need some estimates on size and energy. 

The following Lemmas have been proven in \cite{mtt:fourierbiest}.
 
\begin{lemma}\label{energy-lemma}
Let $j=1,2,3$, $f_j$ be a function in $L^2(\R)$, and let $\Pv$ be a finite collection of tri-tiles.  Then we have
\be{energy-lemma-est}
\modenergy_j((\langle f_j, \phi_{P_j} \rangle)_{P \in \P} ) \lesssim
\| f_j \|_2.
\end{equation}
\end{lemma}

\begin{lemma}\label{size-lemma}
Let $j=1,2,3$, $E_j$ be a set of finite measure, $f_j$ be a function in $X(E_j)$, and let $\Pv$ be a finite collection of tri-tiles.  Then we have
\be{size-lemma-est}
\size_j( (\langle f_j, \phi_{P_j} \rangle)_{\pv \in \Pv} ) \lesssim
\sup_{\pv \in \Pv} \frac{\int_{E_j} \tilde \chi_{I_\pv}^{M}}{|I_\pv|}
\end{equation}
for all $M$, with the implicit constant depending on $M$.
\end{lemma}

We shall also frequently estimate expressions of the form

\begin{equation}\label{bilinear}
| \sum_{P\in\P}a_{P_1} b_{P_2}|
\end{equation}
where $\P$ is a collection of bi-tiles, $a_{P_1}$ are complex 
numbers as before, while $b_{P_2}$ are complex numbers of the form

\begin{equation}\label{bp2}
 b_{P_2} = \langle G\chi_{\{x/N(x)\in\omega_{P_2}\}}, 
\phi_{P_{1}} \rangle.
\end{equation}

We need to define now ``sizes'' and ``energies'' for our 
$b_{P_2}$ sequences. This time, they will no longer depend on the
index ``$j$'' as before.

\begin{definition}\label{bsize-def}
Let $\P$ be a finite collection of bi-tiles, let $\overline{\P}$ be the collection of all bi-tiles corresponding to our fixed dyadic grid (so $\P$ is a finite subset of
$\overline{\P}$) and let 
$(b_{P_2})_{P\in\P}$ be a sequence of complex numbers of the 
form considered above. We define the size of this sequence by

\begin{equation}
\size( (b_{P_2})_{P \in \P} ) :=
\sup_{P\in \P}
\sup_{P'\in \overline{\P}: P<^cP'}
\frac{1}{|I_{P'}|}\int_{\R}|G|
\chi_{\{x/N(x)\in\omega_{P'}\}}\tilde{\chi}_{I_{P'}}^C\,dx
\end{equation}
where $C$ is a fixed big constant. 
We will also need the ``easy variant'' of the size, defined by
\begin{equation}
\size_e( (b_{P_2})_{P \in \P} ) :=
\sup_{P\in \P}
\frac{1}{|I_{P}|}\int_{\R}|G|
\chi_{\{x/N(x)\in\omega_{P}\}}\tilde{\chi}_{I_{P}}^C\,dx.
\end{equation}

Also, since we are in the Fourier setting, we define again a modified energy of this
sequence by

\begin{equation}
\modenergy( (b_{P_2})_{P \in \P} ):= 
\sup_{n\in\Z}\sup_{\D} 2^n (\sum_{P'\in\D} |I_{P'}|)
\end{equation}
where $\D$ ranges over all collections of disjoint bi-tiles $P'\in\overline{\P}$ with the property that there exists $P\in\P$ with $P\leq^c P'$ and such that

$$\int_{\R}|G|
\chi_{\{x/N(x)\in\omega_{P'}\}}\tilde{\chi}_{I_{P'}}^C\,dx \geq 2^n |I_{P'}|.
$$
\end{definition}
As before, $\modenergy( (b_{P_2})_{P \in \P} )$ should be understood as a technical substitute for the more natural (see \cite{mtt:walshbicarleson})
``$\energy( (b_{P_2})_{P \in \P} )$'' defined by

$$\energy( (b_{P_2})_{P \in \P} ):= \sup_{\D\subseteq\overline{\P}}
\sum_{P'\in\D}\int_{\R}|G|
\chi_{\{x/N(x)\in\omega_{P'}\}}\tilde{\chi}_{I_{P'}}^C\,dx,$$
where $\D$ ranges over all collections of disjoint bi-tiles $P'$ for which there exists $P\in\P$ with $P\leq^c P'$.

The main combinatorial tool needed to estimate (\ref{bilinear}) is
the analogue of the above Proposition (\ref{abstract}) for $b_{P_2}$ sequences:
\begin{proposition}\label{babstract}
Let $\P$ be a finite collection of bi-tiles, and for each $P\in\P$ 
let $a_{P_1}$ and $b_{P_2}$ be complex numbers as before. Then,

\begin{equation}\label{bilinear-est}
|\sum_{P\in\P}a_{P_1}b_{P_2}|\lesssim
\size_1((a_{P_1})_{P})^{\theta_1}
\size((b_{P_2})_{P})^{\theta_2}
\modenergy_1((a_{P_1})_{P})^{1-\theta_1}
\modenergy((b_{P_2})_{P})^{1-\theta_2}
\end{equation}
for any $0\leq\theta_1<1$, $0<\theta_2\leq 1$ with $\theta_1+2 \theta_2=1$. Moreover, if the bi-tiles of $\P$ are disjoint, then the inequality holds
even if one replaces $\size((b_{P_2})_{P})$ with the smaller quantity $\size_e((b_{P_2})_{P})$.
\end{proposition}
The proof of this Proposition will be presented later on. 
In the meantime, we will take it for granted.
In order to use Proposition \ref{babstract}, we need again estimates 
on sizes and energies. 

\begin{lemma}\label{benergy-lemma}
Let $f\in L^1(\R)$ and $\P$ be a finite collection of bi-tiles. 
Then, one has

\begin{equation}\label{benergy-lemma-est}
\modenergy((
\langle f\chi_{\{x/N(x)\in\omega_{P_2}\}},\phi_{P_1}\rangle)_{P\in\P})
\lesssim \|f\|_1.
\end{equation}
\end{lemma}

\begin{proof}
Fix $n$ and $\D$ so that the suppremum is attained in the definition of $\modenergy$. Then, we can write
$$\modenergy((
\langle f\chi_{\{x/N(x)\in\omega_{P_2}\}},\phi_{P_1}\rangle)_{P\in\P})\sim 2^n (\sum_{P'\in\D} |I_{P'}| )\lesssim
\sum_{P'\in \D}\int_{\R}|f|
\chi_{\{x/N(x)\in\omega_{P'}\}}\tilde{\chi}_{I_{P'}}^C\,dx$$
and this is smaller than $\|f\|_1$ by using Proposition 3.1 in \cite{laceyt1}.
\end{proof}

\begin{lemma}\label{bsize-lemma}
Let $E$ be a set of finite measure and $f\in X(E)$. Then,

\begin{equation}\label{bsize-lemma-est}
\size((
\langle f\chi_{\{x/N(x)\in\omega_{P_2}\}},\phi_{P_1}\rangle)_{P\in\P})
\lesssim \sup_{P\in\P}\sup_{P'\in\overline{\P}: P\leq^c P'}\frac{\int_E \tilde{\chi}_{I_{P'}}^M}{|I_{P'}|}
\end{equation}
and similarly,

\begin{equation}
\size_e((
\langle f\chi_{\{x/N(x)\in\omega_{P_2}\}},\phi_{P_1}\rangle)_{P\in\P})
\lesssim \sup_{P\in\P}\frac{\int_E \tilde{\chi}_{I_P}^M}{|I_P|}.
\end{equation}
for every $M\leq C$.
\end{lemma}
The proof follows directly from definitions.
In the next section we shall show how the above size and energy estimates can be combined with  Proposition \ref{babstract} 
and the interpolation theory of the previous section to obtain 
Theorem \ref{carleson-fourier}.  
To prove the estimates for the form $\Lambda'_{\P,\Qv}$ 
we need some more sophisticated size and energy estimates, which we will pursue after the proof of Theorem \ref{carleson-fourier}.

\section{Proof of Theorem \ref{carleson-fourier}}\label{carleson-fourier-sec}

We now give a proof of Theorem \ref{carleson-fourier}.  We present its proof
here for expository purposes, and also because we shall need 
Theorem \ref{carleson-fourier} to prove the size and energy estimates needed for Theorem \ref{teorema1}.

First, we linearize the Carleson operator as

$$Cf(x)= \int_{\xi < N(x)} \widehat{f}(\xi) e^{2\pi i x\xi} d\xi$$
and then we dualize it into the bilinear form defined by

$$\Lambda_C(f_1, f_2):= \int_{\R} C(f_1)(x) f_2(x) dx.$$
By standard discretization arguments as in the previous sections, we may reduce the study of $\Lambda_C$ to the study of discretized operators of the form

\bas
 \Lambda_{C, \P}(f_1,f_2) &:= \langle C_{\P}(f_1),f_2 \rangle\\
&= \sum_{P \in \P} 
\langle f_1, \phi_{P_1} \rangle
\langle f_2 \chi_{\{x/N(x)\in\omega_{P_2}\}}, \phi_{P_1} \rangle,
\end{align*}
where $\P$ is some finite collection of bi-tiles.

We shall use the notation of Section \ref{interp-sec}, with the obvious modification for bilinear forms as opposed to trilinear forms.

Let us also consider $E_1, E_2$ sets of finite measure and $1<q<2$.
We are going to prove directly that there exists a major subset
$E'_2$ of $E_2$ so that

\begin{equation}\label{1}
|\Lambda_{C,\P}(f_1,f_2)|\lesssim |E_1|^{1/q}|E_2|^{1/q'},
\end{equation}
for every $f_1\in X(E_1)$, $f_2\in X(E'_2)$ and also that there exists
a major subset $E'_1$ of $E_1$ so that

\begin{equation}\label{2}
|\Lambda_{C,\P}(f_1, f_2)|\lesssim |E_2|,
\end{equation}
for every $f_1\in X(E'_1)$, $f_2\in X(E_2)$. (The reader will notice that a similar argument can be used to prove that there exists a major set $E'_1$ of $E_1$
so that $|\Lambda_{C,\P}(f_1, f_2)|\lesssim \|f_2\|_1$, for every $f_1\in X(E'_1)$ and $f_2\in L^1(\R)$. This is actually equivalent to the fact that the
adjoint $C^*$ of the linearized Carleson operator is of weak type $(1,1)$; thus, it differs at this endpoint from the Carleson operator, which is not of weak type $(1,1)$,
see \cite{fefferman}).

Then, by using the interpolation arguments in \cite{mtt:multilinear}, it follows that
the form $\Lambda_{C,\P}$ is of restricted type $\alpha$, for every $\alpha$
in the interior of the segment defined by the endpoints $(0,1)$ and $(1,0)$.
Finally, Theorem \ref{carleson-fourier} is implied by the classical Marcinkiewicz interpolation
theorem.

It thus remains to prove (\ref{1}) and (\ref{2}).

To prove (\ref{1}), we may assume by scaling invariance, that $|E_2|=1$.
Define the exceptional set

$$\Omega:=\bigcup_{j=1}^2\{x/M\chi_{E_j}>C|E_j|\},$$
where $M$ is the Hardy-Littlewood maximal function. By the classical
Hardy-Littlewood inequality, we have $|\Omega|<1/2$ if $C$ is big enough.
Thus, if we set $E'_2:=E_2\setminus\Omega$, then $E'_2$ is a major subset
of $E_2$. Let now $f_1\in X(E_1)$ and $f_2\in X(E'_2)$. We need to show that

\begin{equation}\label{ab1}
|\sum_{P\in\P}a_{P_1}b_{P_2}|\lesssim |E_1|^{1/q},
\end{equation}
where we denoted by

\begin{equation}
\begin{split}
a_{P_1} &:= \langle f_1,\phi_{P_1} \rangle\\
b_{P_2} &:= \langle f_2\chi_{\{x/N(x)\in\omega_{P_2}\}} ,\phi_{P_1} \rangle.
\end{split}
\end{equation}

We shall make the assumption that

$$1+\frac{\dist(I_P,\Omega^c)}{|I_P|}\sim 2^k$$
for all $P\in \P$, for some $k\geq 0$ independent of $P$ and prove (\ref{1}) with an additional factor of $2^{-\epsilon k}$ on the right hand side (for some $\epsilon >0$).
If we can prove (\ref{1}) in this special case with the indicated gain, then the general case in (\ref{1}) follows by summing in $k$.
 
By the definition of $\Omega$ we thus have

$$ \frac{\int_{E_1}\tilde{\chi}_{I_P}}{|I_P|} \lesssim 2^k |E_1|$$
for all remaining bi-tiles $P \in \P$.  From Lemma \ref{size-lemma} we thus have

$$ \size_1( (a_{P_1})_{P \in \P} ) \lesssim 2^k |E_1|.$$
Also, from Lemma \ref{energy-lemma} and the fact that $f_1 \in X(E_1)$
 we have
$$ \modenergy_1( (a_{P_1})_{P \in \P} ) \lesssim |E_1|^{1/2}.$$
Similarly, for $0\leq k\leq 5$  by usind the definition of $\Omega$ and by applying Lemma \ref{bsize-lemma} we have

$$\size( (b_{P_2})_{P \in \P} )\lesssim 1.$$
For $k>5$, we observe that the corresponding bi-tiles are essentially disjoint and the same Lemma \ref{bsize-lemma}
gives the estimate

$$\size_e( (b_{P_2})_{P \in \P} )\lesssim 2^{-Mk},$$
for any $M>0$. On the other hand, we have from Lemma \ref{benergy-lemma}

$$\modenergy( (b_{P_2})_{P \in \P} )\lesssim 1.$$
From Proposition \ref{babstract} we thus have

$$
|\sum_{P \in \P} 
a_{P_1} b_{P_2}| \lesssim 
(2^k|E_1|)^{\theta_1} (2^{-Mk})^{\theta_2} |E_1|^{(1-\theta_1)/2}= 2^{-k(M\theta_2 - \theta_1)} |E_1|^{(1+\theta_1)/2},
$$
for every $\theta_1\in (0,1)$. If we chose now $\theta_1$ so that
$(1+\theta_1)/2=1/q$ and $M$ big enough, this proves (\ref{1}).

To prove (\ref{2}), we assume again without loss of generality that
$|E_1|=1$, and define $E'_1$ similarly. Then, as before, we restrict the summation over those $P$ having the property that

$$1+\frac{\dist(I_P,\Omega^c)}{|I_P|}\sim 2^k$$
This time, we get the bounds

\begin{equation}
\begin{split}
\size_1( (a_{P_1})_{P \in \P} )& \lesssim 2^{-Mk}\\
\modenergy_1( (a_{P_1})_{P \in \P} )& \lesssim 1\\
\size( (b_{P_2})_{P \in \P} )& \lesssim 2^k |E_2|\\
\energy( (b_{P_2})_{P \in \P} )& \lesssim |E_2|
\end{split}
\end{equation}
and finally, by applying Proposition \ref{babstract}, we obtain

$$|\sum_{P\in\P}a_{P_1} b_{P_2}|\lesssim (2^{-Mk})^{\theta_1} (2^k |E_2|)^{\theta_2}|E_2|^{1-\theta_2}
= 2^{-k(M\theta_1-\theta_2)}|E_2|$$
which again completes the proof, if $M$ is a big constant.

\section{Estimates for $\Lambda'_{P,\Qv}$}

We now continue the study of the form $\Lambda'_{P,\Qv}$.  Fix $P$, $\Qv$
and drop any indices $P$ and $\Qv$ for notational convenience. 

We also drop the $\widetilde{\,\,\,\,\,\,\,\,}$ ' s in the definition of $\Lambda'_{P,\Qv}$, for simplicity.

In the expression $\Lambda'$ the $Q$ tile in the inner summation has a narrower frequency interval, and hence a wider spatial interval, than the $P$ tile in the outer summation.  Thus the inner summation has a poorer spatial localization than the outer sum.  It shall be convenient to reverse the order of summation so that the inner summation is instead more strongly localized spatially than the outer summation.  Specifically, we rewrite $\Lambda'$ as 
$$
\Lambda'_{walsh}(f_1,f_2,f_3)=
\sum_{\vec{Q}\in\Qv}\frac{1}{|I_{\vec{Q}}|^{1/2}} a^{(1)}_{Q_1} a^{(2)}_{Q_2} a^{(3)}_{Q_3}
$$
where
\begin{equation}\label{l'rightform}
\begin{split}
a^{(1)}_{Q_1} &:= \langle f_1,\phi_{Q_1} \rangle\\
a^{(2)}_{Q_2} &:= \langle f_2,\phi_{Q_2} \rangle\\
a^{(3)}_{Q_3} &:= \sum_{P\in\P\,;\,\omega_{Q_3}\subseteq \omega_{P_1}}
\langle f_3\chi_{\{x/N(x)\in\omega_{P_2}\}},\phi_{P_1} \rangle
\langle \phi_{P_1},\phi_{Q_3}\rangle. 
\end{split}
\end{equation}

If $\P'$ is an arbitrary subset of bi-tiles, we also define $C^*_{\P'}$ to be the operator given by

\begin{equation}
C^*_{\P'}(f):= \sum_{P\in\P'}\langle f\chi_{\{x/N(x)\in\omega_{P_2}\}},\phi_{P_1} \rangle
\phi_{P_1}.
\end{equation}
This operator is the adjoint of the Carleson operator and is therefore bounded on every $L^p$ space for $1<p<\infty$ (see Theorem \ref{carleson-fourier}).

To estimate our form $\Lambda'$,
 we need analogues of Lemma \ref{energy-lemma} and Lemma \ref{size-lemma} for $a^{(3)}_{Q_3}$.  One crucial new ingredient in doing so shall be the following 
simple geometric lemma which allows us to decouple the $P$ and $Q$ variables.

\begin{lemma}\label{biest-trick}
Let $i\neq 3$ and let $T\subseteq \Qv$ be an $i$-tree of tri-tiles. For any $\vec{Q}\in T$ we denote by $\P_{\vec{Q}}$ the set

$$\P_{\vec{Q}}:= \{ P\in\P : \omega_{Q_3}\subseteq\omega_{P_1} \}.$$
Similarly, we define the larger set $\P_T$ by

$$\P_T:= \{ P\in\P : \omega_{Q_3}\subseteq\omega_{P_1}, for\,\, some \,\,\vec{Q}\in T  \}.$$
Then, for any $\vec{Q}\in T$ and any function $f$, we have the equality:

\begin{equation}\label{bt}
\langle C^*_{\P_{\vec{Q}}}(f), \phi_{Q_3}\rangle = \langle C^*_{\P_T}(f), \phi_{Q_3}\rangle.
\end{equation}

\end{lemma}

\begin{proof}
Fix $\vec{Q}\in\Qv$ and note that the left hand side of (\ref{bt}) equals

$$\sum_{P\in\P_{\vec{Q}}}\langle f\chi_{\{x/N(x)\in\omega_{P_2}\}},\phi_{P_1} \rangle
\langle \phi_{P_1}, \phi_{Q_3}\rangle,$$
while the right hand side of (\ref{bt}) equals

$$\sum_{P\in\P_T}\langle f\chi_{\{x/N(x)\in\omega_{P_2}\}},\phi_{P_1} \rangle
\langle \phi_{P_1}, \phi_{Q_3}\rangle.$$
Clearly, the sum on the right hand side contains more terms than the sum on the left hand side. Let us now take a look at one
potential non-zero term on the right hand side which does not appear on the left hand side (our claim is that such terms do not exist !).

It corresponds to a bi-tile $P\in\P_T$ so that $\langle\phi_{P_1}, \phi_{Q_3}\rangle \neq 0$. By Plancherel, it follows that
$\frac{9}{10}\omega_{P_1}\cap\frac{9}{10}\omega_{Q_3}\neq\emptyset$ and in particular this implies that $|\omega_{P_1} |\leq 10|\omega_{Q_3} |$
(if not, than we would have $10|\omega_{Q_3} |< |\omega_{P_1} |$ and so $\omega_{Q_3}\subseteq\omega_{P_1}$, which contradicts that the corresponding term
does not appear on the left hand side).

Since this bi-tile $P$ belongs to $\P_T$, it follows that there exists $\vec{Q}'\in T$ such that $\omega_{Q'_3}\subseteq\omega_{P_1}$
and this means that $\omega_{P_1}$ intersects both $\omega_{Q'_3}$ and $\omega_{Q_3}$ and also that

$$| \omega_{Q'_3}|\leq |\omega_{P_1} |\leq 10 |\omega_{Q_3} |.$$
But these two facts contradict that $T$ is an $i$-tree ($i\neq 3$) and the sparseness of our collection of tri-tiles $\Qv$.
The claim follows.

\end{proof}

The energy estimate is given by the following lemma:

\begin{lemma}\label{carleson-energy}
Let $E_3$ be a set of finite measure and $f_3$ be a function in $X(E_3)$.
  Then we have
\be{energy-carleson-est}
\modenergy_3((a^{(3)}_{Q_3})_{\vec{Q} \in \Qv}) \lesssim
|E_3|^{1/2}.
\end{equation}
\end{lemma}

\begin{proof}
Let $n$, $\T$ be an extremizer in the Definition \ref{energy-def} of $\modenergy_3((a^{(3)}_{Q_3})_{\vec{Q} \in \Qv})$. 
By duality, there exists a sequence $(c_{Q_3})_{\vec{Q}}$ of complex numbers
so that

$$\modenergy_3((a^{(3)}_{Q_3})_{\vec{Q} \in \Qv}) \sim \sum_{\vec{Q}\in\T}a^{(3)}_{Q_3} c_{Q_3}.$$
Moreover (see for instance Lemma 6.3 in \cite{mtt:fourierbiest}), the sequence $(c_{Q_3})_{\vec{Q}}$ has the property that

\be{ask-1}
\sum_{\qv \in \tilde T} |c_{Q_3}|^2 \lesssim \frac{|I_{\tilde T}|}{\sum_{T \in \T} |I_T|}
\end{equation}
for all $\tilde T \subseteq T \in \T$. Then, the new expression of the energy becomes

$$\sum_{\vec{Q}\in\T}
 \sum_{P\in\P\,;\,\omega_{Q_3}\subseteq \omega_{P_1}} \langle f_3\chi_{\{x/N(x)\in\omega_{P_2}\}},\phi_{P_1} \rangle
\langle \phi_{P_1},\phi_{Q_3}\rangle c_{Q_3}$$

$$=\sum_{P\in\P}
\langle f_3\chi_{\{x/N(x)\in\omega_{P_2}\}},\phi_{P_1} \rangle 
\langle \phi_{P_1}, \sum_{\vec{Q}\,;\,\omega_{Q_3}\subseteq \omega_{P_1}} c_{Q_3} \phi_{Q_3} \rangle$$

$$:=\sum_{P\in\P} b_{P_2} a_{P_1}.$$
By using Lemma \ref{babstract} in the particular case $\theta_1=0$, $\theta_2=1/2$ the above expression can be majorized by

\begin{equation}\label{gata}
\size((b_{P_2})_{P})^{1/2}
\modenergy_1((a_{P_1})_{P})
\energy((b_{P_2})_{P})^{1/2}.
\end{equation}
Using the previous Lemmas and also Corollary 8.4 in \cite{mtt:fourierbiest}, we obtain the estimates

\begin{equation}
\begin{split}
\modenergy_1( (a_{P_1})_{P \in \P} )& \lesssim 1\\
\size( (b_{P_2})_{P \in \P} )& \lesssim 1\\
\energy( (b_{P_2})_{P \in \P} )& \lesssim |E_3|.
\end{split}
\end{equation}
Using them in (\ref{gata}) we end up with $|E_3|^{1/2}$ as desired.

\end{proof}

We shall also need the following estimate. It is a local version of the Carleson theorem, Theorem \ref{carleson-fourier}.

\begin{lemma}\label{restriction}
Let $i\neq 3$ and let $T\subseteq \Qv$ be an $i$-tree of tri-tiles. Let also $\epsilon >0$. Using the same notations as in Lemma \ref{biest-trick},
the following inequality holds.

$$\left(\int_{\R}|C^*_{\P_T}(f)|^{1+\epsilon}\,\tilde{\chi}_{I_T}^M dx\right)^{1/(1+\epsilon)}\lesssim \left(\int_{E}\tilde{\chi}_{I_T} dx\right)^{1/(1+\epsilon)},$$
for any $f\in X(E)$ and any big $M$.
\end{lemma}

\begin{proof}
To prove our inequality, it is easy to see that it is enough to show that for every interval $I\subseteq \R$ with $|I| = |I_T|$ one has

\begin{equation}\label{111}
\left\|C^*_{\P_T}(f)\right\|_{L^{1+\epsilon}(I)}\lesssim \left(\int_{E}\tilde{\chi}_{I_T} dx\right)^{1/(1+\epsilon)}.
\end{equation}
Then, we estimate the left hand side of (\ref{111}) by

$$\left\|
\sum_{I_P\cap (3I)^c \neq \emptyset}
\langle f\chi_{\{x/N(x)\in\omega_{P_2}\}},\phi_{P_1} \rangle \phi_{P_1}
\right\|_{L^{1+\epsilon}(I)} $$

$$+\left\|
\sum_{I_P\cap (3I)^c = \emptyset}
\langle f\chi_{\{x/N(x)\in\omega_{P_2}\}},\phi_{P_1} \rangle \phi_{P_1}
\right\|_{L^{1+\epsilon}(I)}:= A+B.$$ 
We estimate term A by

$$\sum_{k=0}^{\infty}
\left\|
\sum_{I_P\cap (3I)^c \neq \emptyset ; |I_P|\sim 2^{-k}|I_T|}
\langle f\chi_{\{x/N(x)\in\omega_{P_2}\}},\phi_{P_1} \rangle \phi_{P_1}
\right\|_{L^{1+\epsilon}(I)}$$

$$\lesssim
\sum_{k=0}^{\infty}
\left\|
\sum_{I_P\cap (3I)^c \neq \emptyset ; |I_P|\sim 2^{-k}|I_T|}
\langle |f|, \tilde{\chi}_{I_P}\rangle
\frac{\tilde{\chi}_{I_P}}{|I_P|}
\right\|_{L^{1+\epsilon}(I)}$$

$$\lesssim
\sum_{k=0}^{\infty}
\left(
\sum_{I_P\cap (3I)^c \neq \emptyset ; |I_P|\sim 2^{-k}|I_T|}
(\int_E \tilde{\chi}_{I_P}) (\frac{\dist (I_P, I)}{|I_P|})^{-m}
\right)^{1/(1+\epsilon)}$$

$$\lesssim
\sum_{k=0}^{\infty} 2^{-k}\left(\int_{E}\tilde{\chi}_{I} dx\right)^{1/(1+\epsilon)}\lesssim \left(\int_{E}\tilde{\chi}_{I} dx\right)^{1/(1+\epsilon)}.$$
To estimate term B we have two cases. First, assume that $\supp (f)\subseteq (5I)^c$.
Then we can estimate it by

$$\lesssim
\sum_{k=0}^{\infty}
\left\|
\sum_{I_P\cap (3I)^c = \emptyset ; |I_P|\sim 2^{-k}|I_T|}
\langle |f|, \tilde{\chi}_{I_P}\rangle
\frac{\tilde{\chi}_{I_P}}{|I_P|}
\right\|_{L^{1+\epsilon}(I)}$$

$$\lesssim
\sum_{k=0}^{\infty}
\left(
\sum_{I_P\cap (3I)^c = \emptyset ; |I_P|\sim 2^{-k}|I_T|}
(\int_E \tilde{\chi}_{I_P}) (\frac{\dist (I_P, (5I)^c)}{|I_P|})^{-m}
\right)^{1/(1+\epsilon)}$$

$$\lesssim
\sum_{k=0}^{\infty} 2^{-k}\left(\int_{E}\tilde{\chi}_{I} dx\right)^{1/(1+\epsilon)}\lesssim \left(\int_{E}\tilde{\chi}_{I} dx\right)^{1/(1+\epsilon)}.$$
If, on the other hand, $\supp (f)\subseteq 5I$ then the corresponding term B is smaller
than

$$\left\|
f\right\|_{1+\epsilon} \lesssim \left(\int_{E}\tilde{\chi}_{I} dx\right)^{1/(1+\epsilon)},$$
just by using the fact that the operator $C^*_{\P_T}$ maps $L^{1+\epsilon}$ into itself.

\end{proof}

The analogue of Lemma \ref{size-lemma} is

\begin{lemma}\label{carleson-size}Let $\epsilon>0$,
 $E_3$ be a set of finite measure and $f_3$ be a function in $X(E_3)$.  
Then we have
\be{size-carleson-est}
\size_3((a^{(3)}_{Q_3})_{\vec{Q} \in \Qv}) \lesssim
\sup_{\vec{Q} \in \Qv} (\frac{1}{|I_Q|}\int_{E_3}\tilde{\chi}_{I_Q})^{1/(1+\epsilon)}
\end{equation}
\end{lemma}

\begin{proof}
By Lemma \ref{jn} it suffices to show that
\begin{equation}\label{unuinf}
\frac{1}{|I_T|}
\| (\sum_{\vec{Q} \in T} |a^{(3)}_{Q_3}|^2 
\frac{\chi_{I_Q}}{|I_Q|})^{1/2} \|_{L^{1, \infty}(I_T)}
\lesssim 
\left(\frac{1}{|I_T|}\int_{E_3}\tilde{\chi}_{I_T}\right)^{1/(1+\epsilon)}.
\end{equation}
for some $i \neq 3$ and some $i$-tree $T$. Using Lemma \ref{biest-trick}, the left hand side of (\ref{unuinf})
equals

$$\frac{1}{|I_T|}
\| (\sum_{\vec{Q} \in T} | \langle C^*_{\P_T}(f_3), \phi_{Q_3}\rangle |^2 
\frac{\chi_{I_Q}}{|I_Q|})^{1/2} \|_{L^{1, \infty}(I_T)}.$$
By using Lemma \ref{size-lemma} and Lemma \ref{restriction} this expression can be majorized by

$$\frac{1}{|I_T|}\int_{\R}
| C^*_{\P_T}(f_3)| \tilde{\chi}_{I_T}^M dx \lesssim
\left(\frac{1}{|I_T|}\int_{\R}
| C^*_{\P_T}(f_3)|^{1+\epsilon}\, \tilde{\chi}_{I_T}^M dx\right)^{1/(1+\epsilon)}
\lesssim \left(\frac{1}{|I_T|}\int_{E_3}\tilde{\chi}_{I_T}\right)^{1/(1+\epsilon)},$$
and this ends the proof.

\end{proof}

\section{proof of theorem \ref{teorema1}}

We can now present the proof of Theorem \ref{teorema1}.
Fix the collections $\P$ and $\Qv$.
We first show that $\Lambda'_{\P,\Qv}$
is of restricted weak type $\alpha$ for all admissible 3-tuples 
$(\alpha_1, \alpha_2, \alpha_3)$ 
 arbitrarily close to $A_2$, $A_3$, so that the bad index is 3.

Fix $\alpha$ as above and let $E_1$, $E_2$, $E_3$ be sets of finite measure.  

By scaling invariance we may assume that $|E_3|=1$. 
We need to find a major subset $E'_3$ of $E_3$ such that
$$
|\Lambda'_{\P,\Qv}(f_1, f_2, f_3)|\lesssim |E|^{\alpha}
$$
for all functions $f_i\in X(E'_i)$, $i=1,2,3$.

Define the exceptional set $\Omega$ by
\[\Omega := \bigcup_{j=1}^{3}\{M\chi_{E_j}>C |E_j|\}\]
where $M$ is the dyadic Hardy-Littlewood maximal  function.
By the classical Hardy-Littlewood inequality, we have $|\Omega|<1/2 $
if $C$ is a sufficiently large constant.  Thus if we set 
$E'_3 := E_3 \setminus \Omega$, then $E'_3$ is a major subset of $E_3$.

Let than $f_i \in X(E'_i)$ for $i=1,2,3$.  We need to show that
\begin{equation}\label{22}
|\sum_{\vec{Q} \in \Qv} \frac{1}{|I_{\vec{Q}}|^{1/2}}
a^{(1)}_{Q_1} a^{(2)}_{Q_2} a^{(3)}_{Q_3}| \lesssim |E|^\alpha
\end{equation}
where $a^{(j)}_{Q_j}$ is defined by \eqref{l'rightform}.

We shall make the assumption that

$$1+\frac{\dist(I_{\vec{Q}},\Omega^c)}{|I_{\vec{Q}}|}\sim 2^k$$
for all $\vec{Q}\in \Qv$, for some $k\geq 0$ independent of $\vec{Q}$ and prove (\ref{22}) with an additional factor of $2^{-\epsilon k}$ on the 
right hand side (for some $\epsilon >0$).
If we can prove (\ref{22}) in this special case with the indicated gain, then the general case in (\ref{22}) follows by summing in $k$.

By the definition of $\Omega$ we thus have

$$ \frac{1}{|I_{\vec{Q}}|}\int_{E_j}\tilde{\chi}_{I_{\vec{Q}}} \lesssim 2^k|E_j|$$
for $j=1,2$.  
From Lemma \ref{size-lemma} we thus have
$$ \size_j( (a^{(j)}_{Q_j})_{\vec{Q} \in \Qv} ) \lesssim 2^k|E_j|$$
for $j=1,2$.  
On the other hand, by Lemma \ref{carleson-size} one has
$$ \size_3( (a^{(3)}_{Q_3})_{\vec{Q} \in \Qv} ) \lesssim 2^{-Mk}$$
since $|E_3|=1$, for any big constant $M>0$.
Also, from Lemma \ref{energy-lemma}, Lemma \ref{carleson-energy}
and the fact that $f_j \in X(E'_j)$ we have
$$ \energy_j( (a^{(j)}_{Q_j})_{\vec{Q} \in \Qv} ) \lesssim |E_j|^{1/2},$$
for $j=1,2,3.$
From Proposition \ref{abstract} we thus have
$$
|\sum_{\vec{Q} \in \Qv} \frac{1}{|I_{\vec{Q}}|^{1/2}}
a^{(1)}_{Q_1} a^{(2)}_{Q_2} a^{(3)}_{Q_3}| \lesssim 
\left(\prod_{j=1}^2 |E_j|^{(1-\theta_j)/2} |E_j|^{\theta_j}\right) 2^{-k(M\theta_3 - \theta_1 -\theta_2)} $$
for any $0 \leq \theta_1, \theta_2< 1$ such that there exists
$0< \theta_3<1$ with $\theta_1 + \theta_2 + \theta_3 = 1$.  
The claim then follows by choosing $\theta_1 := 2\alpha_1-1$, $\theta_2 := 2\alpha_2-1$ and $M$ big enough; note that there exist choices of $\alpha$ arbitrarily close to 
$A_2$ or $A_3$, for which the constraints 
on $\theta_1, \theta_2, \theta_3$ are satisfied.  

To prove the restricted type estimates for $\alpha$ arbitrarily
close to $A_4, A_5, A_6, A_1$, one argues in the same way, by taking
advantage of the fact that $\epsilon$ in Lemma \ref{carleson-size}
can be arbitrarily small.

This concludes the proof of Theorem \ref{teorema1}.

\section{Estimates for $\Lambda''_{\P,\Q}$ }

In this section we begin the study of $\Lambda''_{\P,\Q}$.
As before, fix $\P,\Q$ and drop any indices $\P$ and $\Q$ for
notational convenience. We also drop the $\widetilde{\,\,\,\,\,\,\,\,}$ ' s in the definition of $\Lambda''_{P,\Q}$, for simplicity. 

Also, as in the previous sections, it is more 
convenient to rewrite $\Lambda''$ as

\[\Lambda''(f_1,f_2,f_3)=\]
\begin{equation}\label{Lambda''walsh}
\sum_{Q\in\Q}
\langle f_2, \phi_{Q_1}\rangle
\langle \phi_{Q_1} \chi_{\{x/N(x)\in \omega_{Q_2}\}}
\sum_{\omega_{Q_2}\cap\omega_{P_2}\neq\emptyset ; |\omega_{Q_2}| < |\omega_{P_2}|}
\langle f_1, \phi_{P_1}\rangle\phi_{P_1}
\chi_{\{x/N(x)\in \omega_{P_2}\}}, f_3\rangle.
\end{equation}
Expressions of this type have been considered before (see Section 10), but this time, the presence of the inner sum makes their study
much more delicate.

Let us now
fix $E_1,E_2,E_3$ arbitrary sets of finite measure so that $|E_3|=1$.
Define 
\[\Omega := \bigcup_{j=1}^{2}\{M\chi_{E_j}>C|E_j|\}\]
for a large constant $C$, and set $E'_3 := E_3 \setminus \Omega$. Clearly, $|E'_3|>1/2$ if $C$ is a big enough constant.
Pick now $f_1\in X(E_1)$, $f_2\in X(E_2)$ and $f_3\in X(E'_3)$.

The main combinatorial tool needed to estimate our form $\Lambda''$ is
the analogue of the above Proposition \ref{babstract}.
\begin{proposition}\label{ebabstract}
Let $\epsilon >0$ be a small number, and $f_1, f_2, f_3$ as above.
Let also $\Q$ be a finite collection of bi-tiles. Then,

\begin{equation}\label{ebilinear-est}
|\Lambda''(f_1, f_2, f_3)|\lesssim
\size_1((\langle f_2, \phi_{Q_1}\rangle )_{Q})^{\theta_1}
\size((\langle f_3\chi_{\{x : N(x)\in\omega_Q\}},\phi_{Q_1}\rangle )_{Q})^{\theta_2}\cdot
\end{equation}
$$\cdot\modenergy_1((\langle f_2, \phi_{Q_1}\rangle )_{Q})^{1-\theta_1}
\modenergy((\langle f_3\chi_{\{x : N(x)\in\omega_Q\}},\phi_{Q_1}\rangle )_{Q})^{1- \epsilon-\theta_2}\cdot
$$

$$\cdot\left[
\sup_Q(\frac{\int_{E_1}\tilde{\chi}_{I_Q}}{|I_Q|})^{1-\epsilon} + |E_1|^{\alpha}\right]
$$
for any $0 < \theta_1< 1$,  $0< \theta_2 < 1-\epsilon$  with $\theta_1 + 2\theta_2 = 1-2\epsilon$ and $0< \alpha < 1-\epsilon$. Moreover, if the bi-tiles in $\Q$ are disjoint,
then the above inequality holds even when one replaces $\size((\langle f_3\chi_{\{x : N(x)\in\omega_Q\}},\phi_{Q_1}\rangle )_{Q})$ by the smaller quantity
$\size_e((\langle f_3\chi_{\{x : N(x)\in\omega_Q\}},\phi_{Q_1}\rangle )_{Q})$.

\end{proposition}
The proof of this Proposition will be presented later on. 
In the meantime, we will take it for granted.

In the next section we shall show how the above size and energy estimates can be combined with  Proposition \ref{ebabstract} 
and the interpolation theory, to obtain 
Theorem \ref{teorema2}.

\section{Proof of Theorem \ref{teorema2}}

Let $\beta=(\beta_1,\beta_2,\beta_3)$ be an admissible tuple,
very close to either of the points $M_{12}$ or $M_{56}$ or $M_{34}$ or $A$. 

Let us also
fix $E_1,E_2,E_3$ arbitrary sets of finite measure and assume without loss
of generality that $|E_3|=1$.

As usual, we define 
\[\Omega := \bigcup_{j=1}^{2}\{M\chi_{E_j}>C|E_j|\}\]
for a large constant $C$, and set $E'_3 := E_3 \setminus \Omega$.  
We now fix $f_i \in X(E_i)$ for $i=1,2$ and $f_3\in X(E'_3)$.  Our task is then to show that the following inequality

\begin{equation}\label{inegalitatea2}
|\Lambda''(f_1,f_2,f_3)| 
\lesssim |E_1|^{\beta_1} |E_2|^{\beta_2}
\end{equation}
holds, for any $\beta_1$ and $\beta_2$ arbitrarily inside the interval $(0,1)$. Clearly, this would complete the proof.

As before, we may restrict the collection $\Q$ to those bi-tiles $Q$
for which

$$1+\frac{\dist(I_Q, \Omega^c)}{|I_Q|} \sim 2^k$$
for some $k\geq 0$ independent of $Q$ and prove (\ref{inegalitatea2}) with an additional factor of $2^{-\lambda k}$ on the right hand side
(for some $\lambda>0$). If we can prove (\ref{inegalitatea2}) in this special case with the indicated gain, then the general case in
(\ref{inegalitatea2}) follows by summing in $k$.

This implies that
$$ \frac{\int_{E_j}\tilde{\chi}_{I_Q}}{|I_Q|} \lesssim \min (2^k|E_j|, 1)\lesssim (2^k|E_j|)^s$$
for all these tiles $Q \in \Q$, $0<s<1$ and $j=1,2$.  
We also have

$$ \frac{\int_{E'_3}\tilde{\chi}_{I_Q}}{|I_Q|} \lesssim 2^{-Mk}$$
for any big number $M$. We also observe that for $k>5$ the corresponding bi-tiles are essentially disjoint. Using all of these and from Lemma \ref{size-lemma}, 
Lemma \ref{energy-lemma}, we thus have
\bas
\size_1((\langle f_2, \phi_{Q_1}\rangle )_{Q} ) &\lesssim (2^{k}|E_2|)^s\\
\modenergy_1((\langle f_2, \phi_{Q_1}\rangle )_{Q} ) &\lesssim |E_2|^{1/2}\\
\size((\langle f_3\chi_{\{x : N(x)\in\omega_Q\}},\phi_{Q_1}\rangle )_{Q} ) &\lesssim 1\\
\size_e((\langle f_3\chi_{\{x : N(x)\in\omega_Q\}},\phi_{Q_1}\rangle )_{Q} ) &\lesssim 2^{-Mk}\\
\modenergy((\langle f_3\chi_{\{x : N(x)\in\omega_Q\}},\phi_{Q_1}\rangle )_{Q} ) &\lesssim 1.
\end{align*}

By Proposition \ref{ebabstract} we thus can bound the left-hand side of 
\eqref{inegalitatea2} by

$$(2^k|E_2|)^{s\theta_1} (2^{-Mk})^{\theta_2} |E_2|^{(1-\theta_1)/2}(((2^k|E_1|)^{\alpha} + |E_1|^{\alpha})=
2^{-\lambda k} |E_1|^{\alpha} |E_2|^{s\theta_1 - \theta_1/2 + 1/2},
$$
where $\theta_1, \theta_2, \alpha$ are as in Proposition \ref{ebabstract} and $\lambda$ is a number depending on them.
Clearly, $\alpha$ can be chosen arbitrarily inside the interval $(0,1)$ if the $\epsilon$ in Proposition \ref{ebabstract} is small enough. Similarly,
the exponent $s\theta_1 - \theta_1/2 + 1/2$ can be chosen arbitrarily close to $0$ (if $\theta_1$ is close to $1$ and $s$ is close to $0$) and also
arbitrarily close to $1$ (if $\theta_1$ is close to $1$ and $s$ is close to $1$) and  
this finishes the proof, since $\lambda$ can always be made positive if $M$ is chosen big enough.

\section{The tree estimate}

We first recall the following crucial Lemma which is a particular case of Proposition 3.6 in \cite{mptt}.

\begin{lemma}\label{biparameter}
Let $T$ be a $2$-tree of bi-tiles and let $f,g$ be two arbitrary functions. Then,

$$\sum_{P\in T}\left|
\langle f, \phi_{P_1} \rangle
\langle g, \phi_{P_1} \rangle\right|\lesssim
\size_1((\langle f, \phi_{P_1} \rangle)_P)^{1-\theta_1}
\size_1((\langle g, \phi_{P_1} \rangle)_P)^{1-\theta_2}\cdot$$

$$\cdot
\widetilde{\modenergy}((\langle f, \phi_{P_1} \rangle)_P)^{\theta_1}
\widetilde{\modenergy}((\langle g, \phi_{P_1} \rangle)_P)^{\theta_2}$$
for any $0\leq \theta_1, \theta_2 < 1$, $\theta_1 + \theta_2= 1$ where

$$\widetilde{\modenergy}((\langle f, \phi_{P_1} \rangle)_P):=
\sup_{\D\subseteq T}
\left\|
\sum_{P\in\D}\frac{\langle |f|, \tilde{\chi}_{I_P}^C\rangle}{|I_P|}\chi_{I_P}\right\|_{1,\infty},$$
where $\D$ ranges over all subsets of $T$ so that the intervals $\{ I_P : P\in \D \}$ are disjoint.
\end{lemma}

Fix now the collections of bi-tiles $\P$ and $\Q$. We begin our study of the form $\Lambda''_{\P, \Q}$ by considering
the contribution of a single tree (the reader should recall the definition of the form $\Lambda''_{\P, \Q}$ given in Section 14 ).

\begin{lemma}\label{tree-estimate}
Let $\epsilon>0$ a small number, let $T$ be a tree in $\Q$ and $f_1, f_2, f_3$ be three functions as in Proposition \ref{ebabstract}. Then, the following estimate holds:

\begin{equation}\label{tmain}
\left|\Lambda''_{\P, T}(f_1, f_2, f_3)\right|\lesssim
\end{equation}

$$\lesssim
\size_1((\langle f_2, \phi_{Q_1}\rangle )_{Q\in T})
\size((\langle f_3\chi_{\{x : N(x)\in\omega_{Q_2}\}},\phi_{Q_1}\rangle )_{Q\in T})^{1-\epsilon}|I_T|\cdot
$$

$$\cdot
\left[
\sup_Q(\frac{\int_{E_1}\tilde{\chi}_{I_Q}}{|I_Q|})^{1-\epsilon} + |E_1|^{\alpha}\right],
$$
for any $0< \alpha <1-\epsilon$. Moreover, if the bi-tiles in $\Q$ are disjoint, then the expression 
$$\size((\langle f_3\chi_{\{x : N(x)\in\omega_{Q_2}\}},\phi_{Q_1}\rangle )_{Q\in T})$$
can be replaced by the smaller quantity $$\size_e((\langle f_3\chi_{\{x : N(x)\in\omega_{Q_2}\}},\phi_{Q_1}\rangle )_{Q\in T}).$$
\end{lemma}

\begin{proof}
Let $\cal{J}$ be the collection of all maximal intervals $J$ inside our fixed dyadic grid such that $3J$ does not contain any
$I_Q$ with $Q\in T$. Then, clearly, $\cal{J}$ is a partition of the real line $\R$. The left hand side of (\ref{tmain}) can be written as

$$\left|\Lambda''_{\P, T}(f_1, f_2, f_3)\right|$$

$$=\left|\sum_{Q\in T}
\langle f_2, \phi_{Q_1}\rangle
\langle \phi_{Q_1} \chi_{\{x/N(x)\in \omega_{Q_2}\}}
\sum_{\omega_{Q_2}\cap\omega_{P_2}\neq\emptyset ; |\omega_{Q_2}| < |\omega_{P_2}|}
\langle f_1, \phi_{P_1}\rangle\phi_{P_1}
\chi_{\{x/N(x)\in \omega_{P_2}\}}, f_3\rangle\right|$$

$$: =\left|\sum_{Q\in T}
\langle f_2, \phi_{Q_1}\rangle
\langle \phi_{Q_1} \chi_{\{x/N(x)\in \omega_{Q_2}\}}
C_Q(f_1), f_3\rangle\right|$$

$$=\left|\sum_{J\in\cal{J}}\int_J\left(
\sum_{Q\in T}
\langle f_2, \phi_{Q_1}\rangle
\phi_{Q_1} \chi_{\{x/N(x)\in \omega_{Q_2}\}}
C_Q(f_1) f_3\right)(x) dx \right|$$

$$\lesssim\left|\sum_{J\in\cal{J}}\int_J\left(
\sum_{Q\in T : |I_Q| < |J|}
\langle f_2, \phi_{Q_1}\rangle
\phi_{Q_1} \chi_{\{x/N(x)\in \omega_{Q_2}\}}
C_Q(f_1) f_3\right)(x) dx \right|$$

$$+ \left|\sum_{J\in\cal{J}}\int_J\left(
\sum_{Q\in T :|I_Q| > |J| }
\langle f_2, \phi_{Q_1}\rangle
\phi_{Q_1} \chi_{\{x/N(x)\in \omega_{Q_2}\}}
C_Q(f_1) f_3\right)(x) dx \right|$$

$$:= I + II.$$
We first estimate term I. Fix $J\in\cal{J}$ and $Q\in T$ with $|I_Q| < |J|$. We claim that

\begin{equation}\label{t1}
\left|\int_J\left(
\langle f_2, \phi_{Q_1}\rangle
\phi_{Q_1} \chi_{\{x/N(x)\in \omega_{Q_2}\}}
C_Q(f_1) f_3\right)(x) dx \right|
\end{equation}

$$\lesssim\size_1((\langle f_2, \phi_{Q_1}\rangle )_{Q\in T})
\size((\langle f_3\chi_{\{x : N(x)\in\omega_Q\}},\phi_{Q_1}\rangle )_{Q\in T})^{1-\epsilon}
\cdot
$$

$$\cdot\left[
\sup_Q(\frac{\int_{E_1}\tilde{\chi}_{I_Q}}{|I_Q|})^{1-\epsilon} + |E_1|^{\alpha}\right]\cdot
\left(1 + \frac{\dist( I_Q, J)}{|I_Q|}\right)^{-m} |I_Q|,
$$
for any big constant $m$.
Assume that (\ref{t1}) holds. We also have

$$\sum_{Q : |I_Q|\sim 2^k}  
\left(1 + \frac{\dist( I_Q, J)}{|I_Q|}\right)^{-m} |I_Q|\lesssim
2^k \left(1 + \frac{\dist( I_T, J)}{|I_T|}\right)^{-m} 
$$
and after summing over $k$ with $2^k\lesssim |J|$ and over $J\in\cal{J}$ this gives the bound

$$\sum_{J\in\cal{J}} |J| \left(1 + \frac{\dist( I_T, J)}{|I_T|}\right)^{-m} 
\lesssim |I_T|.
$$
This, together with (\ref{t1}) gives the desired estimate. Thus, it remains to prove (\ref{t1}).

We first observe that in order to estimate (\ref{t1}) it is enough to estimate expressions of the form

\begin{equation}\label{t2}
\left|\int_J
\langle f_2, \phi_{Q_1}\rangle
\phi_{Q_1} \chi_{\{x/N(x)\in \omega_{Q_2}\}}
(\sum_{P\in T'}\langle f_1, \phi_{P_1}\rangle \phi_{P_1}) f_3  dx \right|,
\end{equation}
where $P$ runs inside a $2$-tree $T'$ so that all the $|I_P|$' s are smaller than $|I_Q|$.
Now we split the tiles $P$ in $T'$ as $T' = \bigcup_{d\geq 0} T'_d$ where $T'_d$ contains all the tiles having the property that

$$\left(1 + \frac{\dist( I_P, \Omega^c)}{|I_P|}\right)\sim 2^d.$$
We also decompose the interval $J$ as $J = \bigcup_{i=1}^N J_i$ where the $J_i$' s are disjoint intervals so that $|J_i|= |I_Q|$ for every $i=1, ..., N$.
The expression in (\ref{t2}) can be majorized by

$$
\sum_{d\geq 0}\sum_{i=1}^N
\left|\int_{J_i}
\langle f_2, \phi_{Q_1}\rangle
\phi_{Q_1} \chi_{\{x/N(x)\in \omega_{Q_2}\}}
(\sum_{P\in T'_d}\langle f_1, \phi_{P_1}\rangle \phi_{P_1}) f_3  dx \right|.
$$
Fix now $i=1, ..., N$ and $d\geq 0$ and look at the corresponding expression. We decompose it again into

$$\left|\int_{J_i}
\langle f_2, \phi_{Q_1}\rangle
\phi_{Q_1} \chi_{\{x/N(x)\in \omega_{Q_2}\}}
(\sum_{P\in T'_{d, 1}}\langle f_1, \phi_{P_1}\rangle \phi_{P_1}) f_3  dx \right|$$

$$+\left|\int_{J_i}
\langle f_2, \phi_{Q_1}\rangle
\phi_{Q_1} \chi_{\{x/N(x)\in \omega_{Q_2}\}}
(\sum_{P\in T'_{d, 2}}\langle f_1, \phi_{P_1}\rangle \phi_{P_1}) f_3  dx \right|$$

$$:= A + B$$
where $T'_{d, 1}$ contains those tiles $P$ so that $I_P\subseteq 2 J_i$ while $T'_{d, 2}$ contains those tiles so that
$I_P\nsubseteq 2 J_i$.

We now concentrate on $B$. It can be majorized by

$$\size_1((\langle f_2, \phi_{Q_1}\rangle )_{Q\in T})
(\sup_P \frac{\int_{E_1}\tilde{\chi}_{I_P}}{|I_P|} )
\int_{\R}
(\sum_{P\in T'_{d, 2}}\tilde{\chi}_{I_P}^m)
\tilde{\chi}_{I_Q} \chi_{\Omega^c} \chi_{J_i} f_3 dx$$

$$\lesssim
\size_1((\langle f_2, \phi_{Q_1}\rangle )_{Q\in T})
(\sup_P \frac{\int_{E_1}\tilde{\chi}_{I_P}}{|I_P|} )
2^{-dm} 
\left(1 + \frac{\dist( I_Q, J_i)}{|I_Q|}\right)^{-m}\cdot$$

$$\cdot\size((\langle f_3\chi_{\{x : N(x)\in\omega_Q\}},\phi_{Q_1}\rangle )_{Q\in T}) |I_Q|.
$$
Now, after summing over $i= 1, ..., N$ we obtain the bound

$$\size_1((\langle f_2, \phi_{Q_1}\rangle )_{Q\in T})
(\sup_P \frac{\int_{E_1}\tilde{\chi}_{I_P}}{|I_P|} )
2^{-dm}
\left(1 + \frac{\dist( I_Q, J)}{|I_Q|}\right)^{-m}\cdot$$

$$\cdot\size((\langle f_3\chi_{\{x : N(x)\in\omega_Q\}},\phi_{Q_1}\rangle )_{Q\in T}) |I_Q|.
$$

Since we also know that $P\in T'_d$, we have in particular that

$$(\sup_P \frac{\int_{E_1}\tilde{\chi}_{I_P}}{|I_P|} )\lesssim
\min ( 2^d |E_1|, 1)\lesssim 2^{\alpha d} |E_1|^{\alpha},$$
for every $\alpha\in (0,1)$.
Also, since $\size((\langle f_3\chi_{\{x : N(x)\in\omega_Q\}},\phi_{Q_1}\rangle )_{Q\in T})\lesssim 1$ it follows that
$\size((\langle f_3\chi_{\{x : N(x)\in\omega_Q\}},\phi_{Q_1}\rangle )_{Q\in T})
\lesssim
\size((\langle f_3\chi_{\{x : N(x)\in\omega_Q\}},\phi_{Q_1}\rangle )_{Q\in T})^{1-\epsilon}$
and as a consequence, after summing over $d\geq 0$, the new bound is the desired one.

We now concentrate on $A$.
It can be majorized by

\begin{equation}\label{t3}
\size_1((\langle f_2, \phi_{Q_1}\rangle )_{Q\in T})
\sum_{P\in T'_{d,1}}
|\langle f_1, \phi_{P_1}\rangle |
|\langle f_3\tilde{\phi}_{Q_1}\chi_{\{x : N(x)\in\omega_{Q_2}\}}\chi_{J_i}\chi_{\Omega^c}, \phi_{P_1}\rangle |
\end{equation}
where $\tilde{\phi}_{Q_1}:= \phi_{Q_1} |I_Q|^{1/2}$ is an $L^{\infty}$ normalized bump addapted to the interval $I_Q$.
We want to apply Lemma \ref{biparameter} to handle this sum. As before, since $P\in T'_d$ we have the estimate

$$\size_1 ((\langle f_1, \phi_{P_1}\rangle)_P)
\lesssim \sup_P \frac{\int_{E_1}\tilde{\chi}_{I_P}}{|I_P|}\lesssim 
\min ( 2^d |E_1|, 1)\lesssim 2^{\alpha d} |E_1|^{\alpha}.$$
Also, 

$$\size_1((\langle f_3\tilde{\phi}_{Q_1}\chi_{\{x : N(x)\in\omega_{Q_2}\}}\chi_{J_i}\chi_{\Omega^c}, \phi_{P_1}\rangle)_P)
\lesssim 2^{-md}
\left(1 + \frac{\dist( I_Q, J_i)}{|I_Q|}\right)^{-m}.$$
To estimate $\widetilde{\modenergy}((\langle f_1, \phi_{P_1}\rangle)_P)$ fix a set $\D$ as in Lemma \ref{biparameter}
so that the supremum is attained. Since all the $I_P$' s are inside $2 J_i$, we can write

$$\widetilde{\modenergy}((\langle f_1, \phi_{P_1}\rangle)_P)\lesssim
\left\|
\sum_{P\in\D}\frac{\langle |f_1|, \tilde{\chi}_{I_P}^C\rangle}{|I_P|}\chi_{I_P}\right\|_{1,\infty}
$$

$$\lesssim
\left\|
\sum_{P\in\D}\frac{\langle |f_1|\tilde{\chi}_{J_i}, \tilde{\chi}_{I_P}^C\rangle}{|I_P|}\chi_{I_P}\right\|_{1,\infty}
\lesssim 
\| M(|f_1|\tilde{\chi}_{J_i})\|_{1,\infty}$$

$$\lesssim
\| |f_1|\tilde{\chi}_{J_i}\|_{1,\infty}= (\frac{\int_{E_1}\tilde{\chi}_{J_i}}{|J_i|}) |I_Q|\lesssim |I_Q|.$$
Finally,

$$\widetilde{\modenergy}((\langle f_3\tilde{\phi}_{Q_1}\chi_{\{x : N(x)\in\omega_{Q_2}\}}\chi_{J_i}\chi_{\Omega^c}, \phi_{P_1}\rangle)_P)
\lesssim \int_{\R}\tilde{\chi}_{I_Q}\chi_{\{x : N(x)\in\omega_{Q}\}} |f_3| dx
$$

$$=\left(\frac{1}{|I_Q|}\int_{\R}\tilde{\chi}_{I_Q}\chi_{\{x : N(x)\in\omega_{Q}\}} |f_3| dx\right) |I_Q|$$

$$\lesssim \size((\langle f_3\chi_{\{x : N(x)\in\omega_Q\}},\phi_{Q_1}\rangle )_{Q\in T}) |I_Q|.$$
By applying Lemma \ref{biparameter} we estimate (\ref{t3}) by

$$\size_1((\langle f_2, \phi_{Q_1}\rangle )_{Q\in T})
(2^{d\alpha})^{1-\theta_1} \left(2^{-md}\left(1 + \frac{\dist( I_Q, J_i)}{|I_Q|}\right)^{-m}\right)^{1-\theta_2}\cdot
$$

$$\cdot
|I_Q|^{\theta_1} (\size((\langle f_3\chi_{\{x : N(x)\in\omega_{Q_2}\}},\phi_{Q_1}\rangle )_{Q\in T}) |I_Q|)^{\theta_2},
$$
for any $0\leq\theta_1, \theta_2 <1$ with $\theta_1 + \theta_2 = 1$. Now if we choose $\theta_2 = 1-\epsilon$, $\theta_1 = \epsilon$ and $m$ big enough,
we obtain again the desired bound after summing over $d\geq 0$ and $i=1, ..., N$. This ends the discussion on term I.

We now estimate term II. First, we observe that the intervals $J\in \cal{J}$ which contribute to the summation have the property that $J\subseteq 3I_T$.
We then split the tree $T$ as $T = T_1 + T_2$ where $T_1$ is a $1$-tree and $T_2$ is a $2$-tree. As a consequence, our term II also splits 
as

$$II = II_1 + II_2.$$
We first discuss term $II_1$.
We also observe that our tree $T'$ of $P$ bi-tiles also splits as $T':= T'_2 + T'_1$ where $T'_1$ is an $1$-tree and $T'_2$ is a $2$-tree and as before this 
implies a further decomposition of $II_1$ as 

$$II_1 = II_{1a} + II_{1b}.$$
We concentrate on $II_{1a}$ first. We can write it as

$$\sum_J\int_J\left(
\sum_{|I_Q|>|J|}
\langle f_2, \phi_{Q_1}\rangle
\phi_{Q_1} \chi_{\{x/N(x)\in \omega_{Q_2}\}}
\sum_{|\omega_{Q_2}| < |\omega_{P_2}| ; |I_P|<|J|}
\langle f_1, \phi_{P_1}\rangle\phi_{P_1}
\chi_{\{x/N(x)\in \omega_{P_2}\}}f_3\right)dx$$

$$+\sum_J\int_J\left(
\sum_{|I_Q|>|J|}
\langle f_2, \phi_{Q_1}\rangle
\phi_{Q_1} \chi_{\{x/N(x)\in \omega_{Q_2}\}}
\sum_{|\omega_{Q_2}| < |\omega_{P_2}| ; |I_P|>|J|}
\langle f_1, \phi_{P_1}\rangle\phi_{P_1}
\chi_{\{x/N(x)\in \omega_{P_2}\}}f_3\right)dx$$

$$=\sum_J\int_J\left(
\sum_{|I_Q|>|J|} \cdots
\sum_{P\in T'_2 ; |\omega_{Q_2}| < |\omega_{P_2}| ; |I_P|<|J| ; I_P\subseteq 2J}\cdots\right) dx$$

$$+\sum_J\int_J\left(
\sum_{|I_Q|>|J|} \cdots
\sum_{P\in T'_2 ; |\omega_{Q_2}| < |\omega_{P_2}| ; |I_P|<|J| ; I_P\nsubseteq 2J}\cdots\right) dx$$

$$+\sum_J\int_J\left(
\sum_{|I_Q|>|J|} \cdots
\sum_{P\in T'_2 ; |\omega_{Q_2}| < |\omega_{P_2}| ; |I_P|>|J| ; I_P\subseteq 5 I_T}\cdots\right) dx$$

$$+\sum_J\int_J\left(
\sum_{|I_Q|>|J|} \cdots
\sum_{P\in T'_2 ; |\omega_{Q_2}| < |\omega_{P_2}| ; |I_P|>|J| ; I_P\nsubseteq 5 I_T}\cdots\right) dx$$

$$:= II'_{1a} + II''_{1a} + II'''_{1a} + II''''_{1a}.$$
We will treat them one by one. We start with $II'_{1a}$. Fix an interval $J$ and look at the corresponding
function under the integral. It is equal to

$$\left(\sum_{|I_Q|>|J|}
\langle f_2, \phi_{Q_1}\rangle
\phi_{Q_1} \chi_{\{x/N(x)\in \omega_{Q_2}\}}
\right)
\left(
\sum_{P\in T'_2 ; |I_P|<|J| ; I_P\subseteq 2J}\langle f_1, \phi_{P_1}\rangle\phi_{P_1}
\chi_{\{x/N(x)\in \omega_{P_2}\}}\right) f_3 \chi_J \chi_{\Omega^c}$$
and this is pointwise smaller than

$$
\size_1((\langle f_2, \phi_{Q_1}\rangle )_{Q\in T})
\left(
\sum_{P\in T'_2 ; |I_P|<|J| ; I_P\subseteq 2J}\langle f_1, \phi_{P_1}\rangle\phi_{P_1}
\chi_{\{x/N(x)\in \omega_{P_2}\}}\right) f_3 \chi_J \chi_{\Omega^c}.$$
On the other hand, let now $J'$ be an interval in the dyadic grid which contains $J$ and $|J'|= 2|J|$. By the maximality of $J$, it follows
that $3J'$ contains an interval $I_Q$, for some $Q\in T_1$. Then, let $Q_J\in \overline{\Q}$ be a tile with $|I_{Q_J}|= |J'|$ and so that
$Q < Q_J < Q_T$. Clearly, the support of the above function is included inside the set $\{ x : N(x) \in \omega_{Q_J} \}$.
As a consequence of these two facts, we can estimate the integral on $J$ by

$$\size_1((\langle f_2, \phi_{Q_1}\rangle )_{Q\in T})\cdot$$

$$\cdot
\sum_{d\geq 0}
\int_{\R}
\left|
\sum_{P\in T'_2\cap \P_d ; |I_P|<|J| ; I_P\subseteq 2J}\langle f_1, \phi_{P_1}\rangle\phi_{P_1}
\chi_{\{x/N(x)\in \omega_{P_2}\}}
\right|
\chi_{\{ x\in J : N(x) \in \omega_{Q_J} \}}
f_3 \chi_{\Omega^c} dx
$$

\begin{equation}\label{t4}
\lesssim\size_1((\langle f_2, \phi_{Q_1}\rangle )_{Q\in T})\cdot
\end{equation}

$$\cdot
\sum_{d\geq 0}
\int_{\R}
\left|
\sum_{P\in T'_2\cap \P_d ; |I_P|<|J| ; I_P\subseteq 2J}\langle f_1, \phi_{P_1}\rangle\phi_{P_1}
\chi_{\{x/N(x)\in \omega_{P_2}\}}
\right|
\chi_{\{ x : N(x) \in \omega_{Q_J} \}}
f_3 \chi_{\Omega^c} \tilde{\chi}^C_{Q_J}dx.
$$
Fix now $d\geq 0$. To estimate the above integral, it is clearly enough to estimate expressions of the form

$$\int_{\R}
\left(
\sum_{P\in T'_2\cap \P_d ; |I_P|<|J| ; I_P\subseteq 2J}\langle f_1, \phi_{P_1}\rangle\phi_{P_1}
\chi_{\{x/N(x)\in \omega_{P_2}\}}
\right)
\chi_{\{ x : N(x) \in \omega_{Q_J} \}}
f_3 \chi_{\Omega^c} \tilde{\chi}^C_{Q_J} h dx,
$$
where $h\in L^{\infty}$, $\|h\|_{\infty}\leq 1$. This can be further majorized by

\begin{equation}\label{t5}
\sum_{P\in T'_2\cap \P_d ; |I_P|<|J| ; I_P\subseteq 2J}
|\langle f_1, \phi_{P_1}\rangle|
|\langle f_3 \chi_{\Omega^c} \tilde{\chi}^C_{Q_J} h\chi_{\{ x : N(x) \in \omega_{Q_J} \}}, \phi_{P_1}\rangle |.
\end{equation}
To estimate this last expression, we need to apply again Lemma \ref{biparameter}, in the same manner as we did when we estimated (\ref{t3}).
Thus, (\ref{t5}) can be majorized by

$$2^{-\lambda d}
\size((\langle f_3\chi_{\{x : N(x)\in\omega_Q\}},\phi_{Q_1}\rangle )_{Q\in T})^{1-\epsilon} |E_1|^{\alpha} |I_{Q_J}|$$

$$=2^{-\lambda d}
\size((\langle f_3\chi_{\{x : N(x)\in\omega_Q\}},\phi_{Q_1}\rangle )_{Q\in T})^{1-\epsilon} |E_1|^{\alpha} |J|,
$$
where $\lambda$ is a positive number and $0<\alpha <1-\epsilon$. Using this in (\ref{t4}), after summing over $d$ and $J$, we obtain the desired bound.

We now estimate $II''_{1a}$. Just by taking advantage of the decay coming from products of type `` $\phi_{Q_1}\cdot\chi_J$ '', we can easily estimate it
by

$$\sum_{d\geq 0}
2^{-\lambda d}
\sum_J
\size_1((\langle f_2, \phi_{Q_1}\rangle )_{Q\in T})
\size((\langle f_3\chi_{\{x : N(x)\in\omega_Q\}},\phi_{Q_1}\rangle )_{Q\in T})^{1-\epsilon} |E_1|^{\alpha}|Q_J|
$$

$$\lesssim
\sum_J
\size_1((\langle f_2, \phi_{Q_1}\rangle )_{Q\in T})
\size((\langle f_3\chi_{\{x : N(x)\in\omega_Q\}},\phi_{Q_1}\rangle )_{Q\in T})^{1-\epsilon} |E_1|^{\alpha}|J|
$$

$$\lesssim
\size_1((\langle f_2, \phi_{Q_1}\rangle )_{Q\in T})
\size((\langle f_3\chi_{\{x : N(x)\in\omega_Q\}},\phi_{Q_1}\rangle )_{Q\in T})^{1-\epsilon} |E_1|^{\alpha}
|I_T|.
$$
To estimate $II'''_{1a}$ fix again an interval $J$ and look at the corresponding term under the integral. It is given by

$$\sum_{|I_Q|>|J|}
\langle f_2, \phi_{Q_1}\rangle
\phi_{Q_1} \chi_{\{x/N(x)\in \omega_{Q_2}\}}
\sum_{|\omega_{Q_2}| < |\omega_{P_2}| ; |I_P|>|J| ; I_P\subseteq 5 I_T }
\langle f_1, \phi_{P_1}\rangle\phi_{P_1}
\chi_{\{x/N(x)\in \omega_{P_2}\}}f_3$$

$$:= H_J.$$
Fix now $x\in J$. We first observe that since $Q$ ranges inside a tree of type $1$, all the sets of the form
`` $\{x/N(x)\in \omega_{Q_2}\}$ `` are disjoint if the tiles involved have different scales. As a consequence, for our particular fixed $x$,
there is only one $Q$-scale that contributes.
Let us denote this unique scale (which depends on $x$) by $L$. As a consequence, one can write $H_J(x)$ as

\begin{equation}\label{t6}
\left(
\sum_{|I_Q|>|J| ; |\omega_Q| = L}
\langle f_2, \phi_{Q_1}\rangle
\phi_{Q_1}
\right)(x)
\left(
\sum_{|I_P|>|J| ; I_P\subseteq 5 I_T ; L < |\omega_{P_2}|}
\langle f_1, \phi_{P_1}\rangle\phi_{P_1}
\right)(x)
f_3(x).
\end{equation}
Since the intervals $\omega_{P_2}$ are nested, it follows that there is a largest (resp. smallest) interval $\omega_{+}$ (resp. $\omega_{-}$)
of the form $\omega_P$ such that the term in the middle of the product in (\ref{t6}) equals

$$\left(\sum_{ I_P\subseteq 5 I_T ; |\omega_{-}| < |\omega_{P}| < |\omega_{+}|}
\langle f_1, \phi_{P_1}\rangle\phi_{P_1}\right)(x) = 
\left(\sum_{I_P\subseteq 5 I_T }
\langle f_1, \phi_{P_1}\rangle\phi_{P_1}\right)
* \left(\Psi_{+} - \Psi_{-}\right)(x)
$$
where $\Psi_{+}$, $\Psi_{-}$ are well chosen bump functions so that $|\supp (\widehat{\Psi_{+}})|\sim |\omega_{+}|$ and
$|\supp (\widehat{\Psi_{-}})|\sim |\omega_{-}|$. In particular, this implies that this middle term is smaller than

$$\sup_{J\subseteq I}
\frac{1}{|I|}
\int_I
\left|
\sum_{I_P\subseteq 5 I_T }
\langle f_1, \phi_{P_1}\rangle\phi_{P_1}
\right| dy
$$
which is a constant quantity on the interval $J$. As a consequence of thsese observations, our term $II'''_{1a}$
can be estimated by

$$
\size_1((\langle f_2, \phi_{Q_1}\rangle )_{Q})
\sum_J
\left(
\sup_{J\subseteq I}
\frac{1}{|I|}
\int_I
\left|
\sum_{I_P\subseteq 5 I_T }
\langle f_1, \phi_{P_1}\rangle\phi_{P_1}
\right| dy
\right)
\left(\int_J f_3 
\chi_{\{x/N(x)\in \omega_{Q_J}\}} dy\right)
$$

$$\lesssim
\size_1((\langle f_2, \phi_{Q_1}\rangle )_{Q})
\sum_J
\left(
\sup_{J\subseteq I}
\frac{1}{|I|}
\int_I
\left|
\sum_{I_P\subseteq 5 I_T }
\langle f_1, \phi_{P_1}\rangle\phi_{P_1}
\right| dy
\right)
\left(\int_{\R} \tilde{\chi}_{I_{Q_J}}
\chi_{\{x/N(x)\in \omega_{Q_J}\}}f_3  dy\right)
$$

$$\lesssim
\size_1((\langle f_2, \phi_{Q_1}\rangle )_{Q})
\size((\langle f_3\chi_{\{x : N(x)\in\omega_Q\}},\phi_{Q_1}\rangle )_{Q})\cdot
$$

$$\cdot
\sum_J
\left(
\sup_{J\subseteq I}
\frac{1}{|I|}
\int_I
\left|
\sum_{I_P\subseteq 5 I_T }
\langle f_1, \phi_{P_1}\rangle\phi_{P_1}
\right| dy
\right)
|I_{Q_J}|
$$

$$\lesssim
\size_1((\langle f_2, \phi_{Q_1}\rangle )_{Q})
\size((\langle f_3\chi_{\{x : N(x)\in\omega_Q\}},\phi_{Q_1}\rangle )_{Q})\cdot
$$

$$\cdot
\sum_J
\left(
\sup_{J\subseteq I}
\frac{1}{|I|}
\int_I
\left|
\sum_{I_P\subseteq 5 I_T }
\langle f_1, \phi_{P_1}\rangle\phi_{P_1}
\right| dy
\right)
|J|
$$

$$\lesssim
\size_1((\langle f_2, \phi_{Q_1}\rangle )_{Q})
\size((\langle f_3\chi_{\{x : N(x)\in\omega_Q\}},\phi_{Q_1}\rangle )_{Q})\cdot$$

$$\cdot
\left\|
\sum_J
\left(
\sup_{J\subseteq I}
\frac{1}{|I|}
\int_I
\left|
\sum_{I_P\subseteq 5 I_T }
\langle f_1, \phi_{P_1}\rangle\phi_{P_1}
\right| dy
\right)
\chi_{J}
\right\|_{L^1(3I_T)}
$$

$$\lesssim
\size_1((\langle f_2, \phi_{Q_1}\rangle )_{Q})
\size((\langle f_3\chi_{\{x : N(x)\in\omega_Q\}},\phi_{Q_1}\rangle )_{Q})\cdot
$$

$$\cdot\left\|
M(\sum_{I_P\subseteq 5 I_T}
\langle f_1, \phi_{P_1}\rangle\phi_{P_1})
\right\|_{L^1(3I_T)}
$$

$$
\lesssim
\size_1((\langle f_2, \phi_{Q_1}\rangle )_{Q})
\size((\langle f_3\chi_{\{x : N(x)\in\omega_Q\}},\phi_{Q_1}\rangle )_{Q})\cdot
$$

$$\cdot\left\|
M(\sum_{I_P\subseteq 5 I_T}
\langle f_1, \phi_{P_1}\rangle\phi_{P_1})
\right\|_{L^{1/(1-\epsilon)}} |I_T|^{\epsilon}
$$

\begin{equation}\label{t7}
\lesssim
\size_1((\langle f_2, \phi_{Q_1}\rangle )_{Q})
\size((\langle f_3\chi_{\{x : N(x)\in\omega_Q\}},\phi_{Q_1}\rangle )_{Q})\cdot
\end{equation}

$$\cdot\left\|
\sum_{I_P\subseteq 5 I_T}
\langle f_1, \phi_{P_1}\rangle\phi_{P_1}
\right\|_{L^{1/(1-\epsilon)}} |I_T|^{\epsilon}.
$$
Now we also have

$$\left|
\sum_{I_P\subseteq 5 I_T}
\langle f_1, \phi_{P_1}\rangle\phi_{P_1}
\right\|_{L^{1/(1-\epsilon)}} =
\left|
\int_{\R}
\left(
\sum_{I_P\subseteq 5 I_T}
\langle f_1, \phi_{P_1}\rangle\phi_{P_1}
\right)(y) g(y) dy
\right|,
$$
for some $g\in L^{1/\epsilon}$, $\|g\|_{1/\epsilon} =1$. The last term is further smaller than

$$
\sum_{I_P\subseteq 5 I_T}
|\langle f_1, \phi_{P_1}\rangle |
|\langle g, \phi_{P_1}\rangle |
=\int_{\R}
\sum_{I_P\subseteq 5 I_T}
\frac{|\langle f_1, \phi_{P_1}\rangle |}{|I_P|^{1/2}}
\frac{|\langle g, \phi_{P_1}\rangle |}{|I_P|^{1/2}}
\chi_{I_P}(y) dy$$

$$\lesssim
\int_{\R}
\left(
\sum_P
\frac{|\langle f_1, \phi_{P_1}\rangle |^2}{|I_P|}
\chi_{I_P}
\right)^{1/2}
\left(
\sum_P
\frac{|\langle g, \phi_{P_1}\rangle |^2}{|I_P|}
\chi_{I_P}
\right)^{1/2} dy
$$

$$\lesssim
\|
\left(
\sum_P
\frac{|\langle f_1, \phi_{P_1}\rangle |^2}{|I_P|}
\chi_{I_P}
\right)^{1/2}
\|_{L^{1/(1-\epsilon)}}
\cdot
\|
\left(
\sum_P
\frac{|\langle g, \phi_{P_1}\rangle |^2}{|I_P|}
\chi_{I_P}
\right)^{1/2}
\|_{L^{1/\epsilon}}
$$

$$\lesssim
\frac{1}{|I_T|^{1-\epsilon}}
\|
\left(
\sum_P
\frac{|\langle f_1, \phi_{P_1}\rangle |^2}{|I_P|}
\chi_{I_P}
\right)^{1/2}
\|_{L^{1/(1-\epsilon)}}
\cdot 
|I_T|^{1-\epsilon}\lesssim
|I_T|^{1-\epsilon}
(\frac{1}{|I_T|}\int_{E_1}\tilde{\chi}_{I_T})^{1-\epsilon},
$$
by using Lemma \ref{size-lemma}. Inserting this into (\ref{t7}) we obtain the desired bound.

We now estimate $II''''_{1a}$. We first write it (as usual) as

$$\sum_{d\geq 0}\cdots$$
where the $P$-tiles inside `` $\cdots$ `` run inside the set $T'_2\cap\P_d$. Then, fix $d\geq 0$ and $J$ and look at the corresponding
integrand in $II''''_{1a}$.
It is pointwise smaller than

$$\sum_{k\geq 0}
\left|
\sum_{|I_Q|>|J| ; |I_Q|\sim 2^{-k}|I_T|} \cdots
\sum_{P\in T'_2\cap\P_d ; |\omega_{Q_2}| < |\omega_{P_2}| ; |I_P|>|J| ; I_P\nsubseteq 5 I_T}\cdots\right|
\chi_{\Omega^c}\chi_J$$

$$\lesssim2^{-d m}
\size_1((\langle f_2, \phi_{Q_1}\rangle )_{Q})
(\sup_{P\in\P_d}\frac{\int_{E_1}\tilde{\chi}_{I_P}}{|I_P|})
\sum_{k\geq 0}
2^{-k m} f_3 \chi_{\Omega^c}\chi_J\chi_{\{ x : N(x) \in \omega_{Q_J}\} }
$$

$$\lesssim
2^{-d m}
\size_1((\langle f_2, \phi_{Q_1}\rangle )_{Q})
2^{d\alpha} |E_1|^{\alpha} f_3 \chi_J\chi_{\{ x : N(x) \in \omega_{Q_J}\} },
$$
for every $0<\alpha <1$.
Using this information, the integral over $J$ is bounded by (after summing over $d\geq 0$)

$$
\size_1((\langle f_2, \phi_{Q_1}\rangle )_{Q})
|E_1|^{\alpha}
\int_J
\chi_{\{ x : N(x) \in \omega_{Q_J}\} } f_3(x) dx
$$

$$
\lesssim
\size_1((\langle f_2, \phi_{Q_1}\rangle )_{Q})
\size((\langle f_3\chi_{\{x : N(x)\in\omega_{Q_2}\}},\phi_{Q_1}\rangle )_{Q})
|E_1|^{\alpha} |J|
$$
and this, after summing over $J\in\cal{J}$, gives the desired bound.
To finish the discussion on term $II_1$ we need to discuss now term $II_{1b}$. We first split it 
as

$$\sum_{d\geq 0}\cdots$$
as before. Since now both our trees $T$ and $T'$ are $1$-trees, it follows that the sets
$\{x : N(x)\in\omega_{Q_2}\}$ and $\{x : N(x)\in\omega_{P_2}\}$ are disjoint, if they correspond to different scales. As a consequence,
for a fixed $d\geq 0$ and $J\in\cal{J}$, the corresponding integrand is pointwise smaller than

$$
\size_1((\langle f_2, \phi_{Q_1}\rangle )_{Q})
(\sup_{P\in\P_d}\frac{\int_{E_1}\tilde{\chi}_{I_P}}{|I_P|})
2^{-d m}
 f_3
\chi_J\chi_{\{ x : N(x) \in \omega_{Q_J}\} }
$$

$$\lesssim
\size_1((\langle f_2, \phi_{Q_1}\rangle )_{Q})
2^{d\alpha}
|E_1|^{\alpha}
2^{-d m}
 f_3
\chi_J\chi_{\{ x : N(x) \in \omega_{Q_J}\} }
$$
for every $0<\alpha <1$. Integrating this over $J$ and summing over $d\geq 0$ we obtain the bound

$$\size_1((\langle f_2, \phi_{Q_1}\rangle )_{Q})
\size((\langle f_3\chi_{\{x : N(x)\in\omega_{Q_2}\}},\phi_{Q_1}\rangle )_{Q})
|E_1|^{\alpha} |J|
$$
which, after summing over $J\in\cal{J}$, becomes the desired bound.

It remains to estimate term $II_2$. First of all, we write as before 
$T' = T'_2 + T'_1$ where $T'_1$ is an $1$-tree and $T'_2$ is a $2$-tree and this implies a decomposition of $II_2$
as

$$II_2 = II_{2a} + II_{2b}.$$
Furthermore, we split again the first term $II_{2a}$ as

$$II_{2a} = II'_{2a} + II''_{2a} + II'''_{2a} + II''''_{2a}$$
where each of these terms correspond to the same summation constraints as before, when we decomposed the term $II_{1a}$.

We now estimate $II'_{2a}$. As usual we split it as $\sum_{d\geq 0}\cdots$. Then, we fix $d\geq 0$ and $J$ and look at the corresponding integrand.
It can be written as

$$
\left(
\sum_{|I_Q|>|J|}
\langle f_2, \phi_{Q_1}\rangle
\phi_{Q_1} \chi_{\{x/N(x)\in \omega_{Q_2}\}}\right)
\left(
\sum_{|I_P|<|J| ; I_P\subseteq 2J}
\langle f_1, \phi_{P_1}\rangle\phi_{P_1}
\chi_{\{x/N(x)\in \omega_{P_2}\}}f_3
\chi_{\Omega^c}\chi_J
\right)$$

$$\lesssim
\left(
\sup_{J\subseteq I}
\frac{1}{|I|}
\int_I
\left|
\sum_{Q}
\langle f_2, \phi_{Q_1}\rangle\phi_{Q_1}
\right| dy
\right)\cdot
$$

$$\cdot
\left|
\sum_{|I_P|<|J| ; I_P\subseteq 2J}
\langle f_1, \phi_{P_1}\rangle\phi_{P_1}
\chi_{\{x/N(x)\in \omega_{P_2}\}}f_3
\chi_{\Omega^c}\chi_J\chi_{\{x/N(x)\in \omega_{Q_J}\}}\right|
$$
by using the same geometric arguments used to estimate term $II'''_{2a}$.
Integrating over $J$ this gives the bound

\begin{equation}\label{t8}
\left(
\sup_{J\subseteq I}
\frac{1}{|I|}
\int_I
\left|
\sum_{Q}
\langle f_2, \phi_{Q_1}\rangle\phi_{Q_1}
\right| dy
\right)\cdot
\sum_{|I_P|<|J| ; I_P\subseteq 2J}
|\langle f_1, \phi_{P_1}\rangle |
|\langle f_3 \chi_{\Omega^c} h \chi_J \chi_{\{x/N(x)\in \omega_{Q_J}\}}, \phi_{P_1}\rangle|,
\end{equation}
for some $h\in L^{\infty}$, $\|h\|_{\infty}\leq 1$.
The last sum can be estimated as before (using Lemma \ref{biparameter}) by

$$(2^{d\alpha} |E_1|^{\alpha})^{1-\epsilon}(2^{-m d})^{\epsilon} |J|^{\epsilon}
(\size((\langle f_3\chi_{\{x : N(x)\in\omega_{Q_2}\}},\phi_{Q_1}\rangle )_{Q})|J|)^{1-\epsilon}
$$

$$= 2^{-\lambda d}
|E_1|^{\alpha(1-\epsilon)}
(\size((\langle f_3\chi_{\{x : N(x)\in\omega_{Q_2}\}},\phi_{Q_1}\rangle )_{Q}))^{1-\epsilon}
|J|,
$$
for some $\lambda >0$, if $m$ is big enough.
As a consequence, after summing over $d\geq 0$ the term $II'_{2a}$
becomes smaller than

$$|E_1|^{\alpha(1-\epsilon)}
(\size((\langle f_3\chi_{\{x : N(x)\in\omega_{Q_2}\}},\phi_{Q_1}\rangle )_{Q}))^{1-\epsilon}
\sum_{J}
\left(
\sup_{J\subseteq I}
\frac{1}{|I|}
\int_I
\left|
\sum_{Q}
\langle f_2, \phi_{Q_1}\rangle\phi_{Q_1}
\right| dy
\right) |J|
$$
and this, by an argument used before is smaller than

$$|E_1|^{\alpha(1-\epsilon)}
(\size((\langle f_3\chi_{\{x : N(x)\in\omega_{Q_2}\}},\phi_{Q_1}\rangle )_{Q}))^{1-\epsilon}
\size_1((\langle f_2, \phi_{Q_1}\rangle )_{Q}) |I_T|.
$$
We now estimate $II''_{2a}$. As before, we first we decompose it as $\sum_{d\geq 0}\cdots$. Then, we fix $d\geq 0$ and $J$ and look at the corresponding
integrand. It can be written as

\begin{equation}
\left(
\sum_{|I_Q|>|J|}
\langle f_2, \phi_{Q_1}\rangle
\phi_{Q_1} \chi_{\{x/N(x)\in \omega_{Q_2}\}}\right)
\left(
\sum_{|I_P|<|J| ; I_P\nsubseteq 2J}
\langle f_1, \phi_{P_1}\rangle\phi_{P_1}
\chi_{\{x/N(x)\in \omega_{P_2}\}}f_3
\chi_{\Omega^c}\chi_J
\right)
\end{equation}

$$\lesssim
\left(
\sup_{J\subseteq I}
\frac{1}{|I|}
\int_I
\left|
\sum_{Q}
\langle f_2, \phi_{Q_1}\rangle\phi_{Q_1}
\right| dy
\right)
(\sup_{P\in\P_d}\frac{\int_{E_1}\tilde{\chi}_{I_P}}{|I_P|})
2^{-m d} f_3 \chi_J\chi_{\{x/N(x)\in \omega_{Q_J}\}}
$$

$$\lesssim
\left(
\sup_{J\subseteq I}
\frac{1}{|I|}
\int_I
\left|
\sum_{Q}
\langle f_2, \phi_{Q_1}\rangle\phi_{Q_1}
\right| dy
\right)
2^{d\alpha} |E_1|^{\alpha} 2^{-m d} f_3 \chi_J
\chi_{\{x/N(x)\in \omega_{Q_J}\}}
$$
for $0<\alpha<1$. After we integrate over $J$ we get the bound

$$
2^{-\lambda d} |E_1|^{\alpha}
\left(
\sup_{J\subseteq I}
\frac{1}{|I|}
\int_I
\left|
\sum_{Q}
\langle f_2, \phi_{Q_1}\rangle\phi_{Q_1}
\right| dy
\right)
\size((\langle f_3\chi_{\{x : N(x)\in\omega_{Q_2}\}},\phi_{Q_1}\rangle )_{Q})|J|.
$$
Summing now over $d$ and $J$ we obtain the desired bound.

We now estimate term $II'''_{2a}$ as follows. Fix $J\in \cal{J}$. The corresponding integrand can be written as

\begin{equation}\label{t9}
\sum_{|I_Q|>|J|}
\langle f_2, \phi_{Q_1}\rangle
\phi_{Q_1} \chi_{\{x/N(x)\in \omega_{Q_2}\}}
\sum_{|\omega_{Q_2}| < |\omega_{P_2}|  |I_P|>|J| ; I_P\subseteq 5I_T}
\langle f_1, \phi_{P_1}\rangle\phi_{P_1}
\chi_{\{x/N(x)\in \omega_{P_2}\}}f_3
\chi_J.
\end{equation}
Using a similar analysis on the geometry of the frequency intervals, as in the case when we estimated term $II'''_{1a}$,
the above term (\ref{t9}) is pointwise smaller than

$$
\sup_{J\subseteq I}
\frac{1}{|I|}
\int_I
\left|
\Pi (\sum_{Q}
\langle f_2, \phi_{Q_1}\rangle\phi_{Q_1}, \sum_P \langle f_1, \phi_{P_1}\rangle\phi_{P_1})
\right| 
$$

$$+ 
\left(
\sup_{J\subseteq I}
\frac{1}{|I|}
\int_I
\left|
\sum_{Q}
\langle f_2, \phi_{Q_1}\rangle\phi_{Q_1}
\right|
\right)
\left(
\sup_{J\subseteq I}
\frac{1}{|I|}
\int_I
\left|
\sum_{Q}
\langle f_1, \phi_{P_1}\rangle\phi_{P_1}
\right|
\right),
$$
where $\Pi$ is a paraproduct well adapted tp the frequency intervals of the trees $T_2$ and $T'_2$. Notice that this expression is constant on the interval $J$.
After integrating over $J$ and summing over all the intervals in $\cal{J}$ we can bound term $II'''_{2a}$ as before by

$$
\size((\langle f_3\chi_{\{x : N(x)\in\omega_{Q_2}\}},\phi_{Q_1}\rangle )_{Q})\cdot
\left\|
M(\Pi (\sum_{Q}
\langle f_2, \phi_{Q_1}\rangle\phi_{Q_1}, \sum_P \langle f_1, \phi_{P_1}\rangle\phi_{P_1}) )
\right\|_{L^1(3I_T)}
$$

$$+ 
\size((\langle f_3\chi_{\{x : N(x)\in\omega_{Q_2}\}},\phi_{Q_1}\rangle )_{Q})\cdot
\left\|
M(\sum_{Q}
\langle f_2, \phi_{Q_1}\rangle\phi_{Q_1} )M(\sum_{P}
\langle f_1, \phi_{P_1}\rangle\phi_{P_1} )
\right\|_{L^1(3I_T)}
$$

$$\lesssim
\size((\langle f_3\chi_{\{x : N(x)\in\omega_{Q_2}\}},\phi_{Q_1}\rangle )_{Q})\cdot
\left\|
\Pi (\sum_{Q}
\langle f_2, \phi_{Q_1}\rangle\phi_{Q_1}, \sum_P \langle f_1, \phi_{P_1}\rangle\phi_{P_1})
\right\|_{L^{1+\epsilon}(3I_T)}|I_T|^{\epsilon/(1+\epsilon)}
$$

$$+
\size((\langle f_3\chi_{\{x : N(x)\in\omega_{Q_2}\}},\phi_{Q_1}\rangle )_{Q})\cdot
\left\|
\sum_{Q}
\langle f_2, \phi_{Q_1}\rangle\phi_{Q_1}
\right\|_{L^{1/\epsilon}}\cdot
\left\|
\sum_{P}
\langle f_1, \phi_{P_1}\rangle\phi_{P_1}
\right\|_{L^{1/(1-\epsilon)}}
$$

$$\lesssim
\size((\langle f_3\chi_{\{x : N(x)\in\omega_{Q_2}\}},\phi_{Q_1}\rangle )_{Q})
|I_T|^{\epsilon/(1+\epsilon)}
\left\|
\sum_{Q}
\langle f_2, \phi_{Q_1}\rangle\phi_{Q_1}
\right\|_{L^{(1+\epsilon)/\epsilon^2}}\cdot
\left\|
\sum_{P}
\langle f_1, \phi_{P_1}\rangle\phi_{P_1}
\right\|_{L^{1/(1-\epsilon)}}
$$

$$+
\size((\langle f_3\chi_{\{x : N(x)\in\omega_{Q_2}\}},\phi_{Q_1}\rangle )_{Q})\cdot
\left\|
\sum_{Q}
\langle f_2, \phi_{Q_1}\rangle\phi_{Q_1}
\right\|_{L^{1/\epsilon}}\cdot
\left\|
\sum_{P}
\langle f_1, \phi_{P_1}\rangle\phi_{P_1}
\right\|_{L^{1/(1-\epsilon)}}
$$

$$\lesssim
\size((\langle f_3\chi_{\{x : N(x)\in\omega_{Q_2}\}},\phi_{Q_1}\rangle )_{Q})
\size_1((\langle f_2, \phi_{Q_1}\rangle )_{Q})
(\frac{1}{|I_T|}\int_{E_1}\tilde{\chi}_{I_T})^{1-\epsilon} |I_T|,
$$ 
which is what we wanted. Term $II''''_{2a}$ can be estimated in a similar way with term $II''''_{1a}$ and is left to the reader.
It remains to estimate term $II_{2b}$ in order to finish the whole proof. We first split it as usual, as $\sum_{d\geq 0}\cdots$.
We then fix $d\geq 0$ and $J$ and look at the corresponding integrand. Since $P$ runs now inside a tree of type $1$, the intervals
$\{x : N(x)\in\omega_{P_2}\}$ are disjoint if they correspond to different scales. In particular, for a fixed $x\in J$
there is only one scale that contributes to our summation. As a consequence, we can pointwise estimate our integrand by

\begin{equation}
\left(
\sup_{J\subseteq I}
\frac{1}{|I|}
\int_I
\left|
\sum_{Q}
\langle f_2, \phi_{Q_1}\rangle\phi_{Q_1}
\right| dy
\right)
(\sup_{P\in\P_d}\frac{\int_{E_1}\tilde{\chi}_{I_P}}{|I_P|})
2^{-m d} f_3 \chi_J\chi_{\{x/N(x)\in \omega_{Q_J}\}}
\end{equation}
and this, as we have already seen, leads to the right estimate. This ends the proof of our inequality.

In the particular case when all the bi-tiles in $\Q$ are disjoint, every tree is an ``one bi-tile tree'' and the whole
proof becomes much simpler. As a consequence, one easily observes that $\size[...]$ can be replaced by $\size_e[...]$.
The proof is now complete.

\end{proof}

\section{Combinatorial lemmas}

In order to prove Propositions \ref{babstract} and \ref{ebabstract} we need to recall certain standard combinatorial Lemmas.
To bootstrap the summation over a single tree $T$ as in Lemma \ref{tree-estimate}, to a summation over the whole $\Q$,
 we would like to partition $\Q$ into trees $T$ for which one has good control over $\sum_T |I_T|$.  
This will be accomplished by several decomposition lemmas. The first
one appeared in \cite{mtt:fourierbiest} (see Proposition 12.2.).

\begin{proposition}\label{decomp}  Let $j = 1,2,3$, 
$\Q'$ be a subset of $\Q$, $n \in \Z$, $f$ be a function and suppose that
$$ \size_j((\langle f,\phi_{Q_j}\rangle)_{Q \in \Q'} ) \leq 2^{-n} 
\modenergy_j((\langle f,\phi_{Q_j}\rangle)_{Q \in \Q} ).$$
Then we may decompose $\Q' = \Q'' \cup \Q'''$ such that
\be{size-lower}
\size_j((\langle f,\phi_{Q_j}\rangle)_{Q \in \Q''} ) \leq 2^{-n-1} 
\modenergy_j((\langle f,\phi_{Q_j}\rangle)_{Q \in \Q} )
\end{equation}
and that $\Q'''$ can be written as the disjoint union of trees $\T$ such that
\be{tree-est}
\sum_{T \in \T} |I_T| \lesssim 2^{2n}.
\end{equation}
\end{proposition}
By iterating this proposition one obtains (see again \cite{mtt:fourierbiest}, Corollary 12.3).

\begin{corollary}\label{decomp-cor}  There exists a partition
$$ \Q = \bigcup_{n \in \Z} \Q_n$$
where for each $n \in \Z$ and $j = 1,2,3$ we have
$$ \size_j((\langle f,\phi_{Q_j}\rangle)_{Q \in \Q_n} ) 
\leq 
\min(2^{-n} \modenergy_j((\langle f,\phi_{Q_j}\rangle)_{Q \in \Q} ) , 
\size_j((\langle f,\phi_{Q_j}\rangle)_{Q \in \Q} )).$$
Also, we may cover $\Q_n$ by a collection $\T_n$ of trees such that
$$ \sum_{T \in \T_n} |I_T| \lesssim 2^{2n}.$$
\end{corollary}
The next Proposition together with its corollary are also known (see \cite{laceyt1}, Proposition 3.1).

\begin{proposition}\label{decompc}
 Let  
$\Q'$ be a subset of $\Q$, $n \in \Z$, $f$ be a function and suppose that
$$ \size((\langle f\chi_{\{x/N(x)\in\omega_{Q_2}\}},\phi_{Q_1}
\rangle)_{Q\in\Q'}) \leq 2^{-n} 
\modenergy((\langle f\chi_{\{x/N(x)\in\omega_{Q_2}\}},\phi_{Q_1}
\rangle)_{Q\in\Q})
.$$
Then we may decompose $\Q' = \Q'' \cup \Q'''$ such that
\be{size-lowerc}
\size((\langle f\chi_{\{x/N(x)\in\omega_{Q_2}\}},\phi_{Q_1}
\rangle)_{Q\in\Q''})
\leq 2^{-n-1}
\modenergy((\langle f\chi_{\{x/N(x)\in\omega_{Q_2}\}},\phi_{Q_1}
\rangle)_{Q\in\Q})
\end{equation}
and that $\Q'''$ can be written as the disjoint union of trees $\T$ such that
\be{tree-estc}
\sum_{T \in \T} |I_T| \lesssim 2^{n}.
\end{equation}
Moreover, if all the bi-tiles in $\Q$ are disjoint, the Proposition holds if one replaces $\size[...]$ by $\size_e[...]$.
\end{proposition}

\begin{corollary}\label{decompc-cor}  There exists a partition
$$ \Q= \bigcup_{n \in \Z} \Q_n$$
where for each $n \in \Z$  we have
$$ \size((\langle f\chi_{\{x/N(x)\in\omega_{Q_2}\}},\phi_{Q_1}
\rangle)_{Q\in\Q_n})
\leq $$
$$
\min(2^{-n} 
\modenergy((\langle f\chi_{\{x/N(x)\in\omega_{Q_2}\}},\phi_{Q_1}
\rangle)_{Q\in\Q}),
\size((\langle f\chi_{\{x/N(x)\in\omega_{Q_2}\}},\phi_{Q_1}
\rangle)_{Q\in\Q}) )
.$$
Also, we may cover $\Q_n$ by a collection $\T_n$ of trees such that
$$ \sum_{T \in \T_n} |I_T| \lesssim 2^{n}.$$
\end{corollary}

\section{Proof of Propositions \ref{babstract} and \ref{ebabstract}}

It remains to present the proofs of Propositions \ref{babstract} and \ref{ebabstract}.
We start with Proposition \ref{babstract}. Fix $\Q$ a collection of bi-tiles and let $a_{Q_1}$ and $b_{Q_2}$ be complex numbers as before, given by

$$a_{Q_1}:= \langle f_1, \phi_{Q_1}\rangle$$
and
$$b_{Q_2}:= \langle f_2 \chi_{\{x/N(x)\in\omega_{Q_2}\}},\phi_{Q_1}\rangle.$$
Fix also $\theta_1, \theta_2 \in (0, 1)$ so that $\theta_1+ 2\theta_2=1$. First of all, let us recall the following
standard estimate (see \cite{laceyt1}):

\begin{equation}\label{oldtree}
\left|
\sum_{Q\in T}a_{Q_1} b_{Q_2}
\right|\lesssim
\size_1((a_{Q_1})_{Q\in T}) \size((b_{Q_2})_{Q\in T}) |I_T|
\end{equation}
for any tree $T$ in $\Q$.
During the proof, we will write for simplicity 
$S_1:= \size_1((a_{Q_1})_{Q\in \Q})$, $S_2:= \size((b_{Q_2})_{Q\in \Q})$,
$E_1:= \modenergy_1((a_{Q_1})_{Q\in \Q} )$ and $E_2:= \modenergy((b_{Q_2})_{Q\in \Q} )$.

If we apply Corollaries \ref{decomp-cor} and \ref{decompc-cor} to the functions $\frac{f_1}{E_1}$ and $\frac{f_2}{E_2}$ respectively,
we obtain decompositions

$$\Q = \bigcup_{n \in \Z} \Q_n^j$$
for $j=1,2$ such that each $\Q_n^j$ can be written as a union of subsets in $\T_n^j$ satisfying the properties of those Corollaries for $j=1,2$.
In particular, we can write the left hand side of our desired inequality as

\begin{equation}\label{*}
E_1 E_2 \sum_{n_1, n_2} \sum_{T\in\T^{n_1, n_2}}\left|
\sum_{Q\in T}a_{Q_1} b_{Q_2}
\right|
\end{equation}
where $\T^{n_1, n_2}:= \T^1_{n_1}\cap \T^2_{n_2}$. By using the above tree estimate (\ref{oldtree}) one can majorize (\ref{*}) by

\begin{equation}\label{**}
E_1 E_2 \sum_{n_1, n_2}2^{-n_1} 2^{-n_2} \sum_{T\in\T^{n_1, n_2}}|I_T|,
\end{equation}
where, according to the same Corollaries the summation goes over those $n_1, n_2\in\Z$ such that

$$2^{-n_j}\lesssim \frac{S_j}{E_j}$$
for $j=1,2$. On the other hand we also know that we can estimate the inner sum in (\ref{**}) in two different ways,
namely

\begin{equation}\label{tn1n2}
\sum_{T\in\T^{n_1, n_2}}|I_T|\lesssim
2^{2n_1},\,\,2^{n_2}
\end{equation}
and so, in particular we can also write

\begin{equation}\label{***}
\sum_{T\in\T^{n_1, n_2}}|I_T|\lesssim 2^{2n_1\alpha_1} 2^{n_2\alpha_2}
\end{equation}
for any $0\leq \alpha_1, \alpha_2\leq 1$ with $\alpha_1 + \alpha_2 = 1$. Using (\ref{***}), we can estimate (\ref{**}) further by

$$E_1 E_2\sum_{n_1, n_2} 2^{-n_1(1-2\alpha_1)} 2^{-n_2(1-\alpha_2)}
$$

$$\lesssim E_1 E_2 (\frac{S_1}{E_1})^{1-2\alpha_1} (\frac{S_2}{E_2})^{1-\alpha_2}=
S_1^{1-2\alpha_1} S_2^{1-\alpha_2} E_1^{2\alpha_1} E_2^{\alpha_2},$$
as long as $1-2\alpha_1>0$, $1-\alpha_2>0$ and $\alpha_1 + \alpha_2 = 1$.
Now, if we choose $\alpha_1, \alpha_2$ such that $\theta_1=1-2\alpha_1$ and $\theta_2=1-\alpha_2$ we observe that $\theta_1+ 2\theta_2=1$ and the last term above
becomes

$$S_1^{\theta_1} S_2^{\theta_2} E_1^{1-\theta_1} E_2^{1-\theta_2}.$$
This proves our inequality in the case when $\theta_1, \theta_2 \in (0, 1)$.
We now prove the endpoint case $\theta_1=0$ (and so $\theta_2=1/2$).

From (\ref{tn1n2}) we have

$$\sum_{T\in\T^{n_1, n_2}}|I_T|\lesssim
\min (2^{2n_1}, 2^{n_2})$$
and as a consequence, our sum in (\ref{**}) can be estimated by

$$E_1 E_2 \sum_{n_1, n_2}2^{-n_1} 2^{-n_2}\min (2^{2n_1}, 2^{n_2})$$

$$=E_1 E_2\sum_{n_2} 2^{-n_2} \sum_{n_1}\min (2^{n_1}, 2^{-n_1}2^{n_2})$$

$$\lesssim E_1 E_2 \sum_{n_2} 2^{-n_2} 2^{\frac{n_2}{2}}= E_1 E_2 \sum_{n_2}2^{-\frac{n_2}{2}}$$

$$= E_1 E_2 (\frac{S_2}{E_2})^{1/2} =  S_2^{1/2} E_1 E_2^{1/2}$$
and this ends the proof of the main part of the Proposition. 

In the particular case when all the bi-tiles in $\Q$ are disjoint, one just has to observe that all the trees $T$ in $\Q$ are ``one bi-tile trees'' and 
then to use the trivial inequality 
\begin{equation}
\left|
\sum_{Q\in T}a_{Q_1} b_{Q_2}
\right|\lesssim
\size_1((a_{Q_1})_{Q\in T}) \size_e((b_{Q_2})_{Q\in T}) |I_T|
\end{equation}
instead of the previous (\ref{oldtree}).

The proof of the remaining Proposition \ref{ebabstract} is very similar and will be omitted. The only difference is that one has to use Lemma \ref{tree-estimate} instead
of the tree estimate (\ref{oldtree}). 
The extra term ``$[...]$'' (which did not appear in the proof of Proposition \ref{babstract} ) is harmless in the whole process and can be factored out. 
This is why it is the same in the 
statements of Lemma \ref{tree-estimate} and Proposition \ref{ebabstract}. Finally, when one keeps track of the ``numerology'', one ends up with a condition
depending on $\epsilon$, instead of the previous $\theta_1+ 2\theta_2=1$.

\end{document}